\documentclass{tran-l}
\usepackage{amssymb}

\usepackage{amsmath}

\input{tcilatex}
\theoremstyle{definition}
\theoremstyle{remark}
\numberwithin{equation}{section}

\input tcilatex

\begin{document}
\title[A residue of complex function...]{A residue of complex function in three-dimensional vector space}
\author{\textbf{branko saric}}
\address{The \textit{Institute ''Kirilo Savic``, 11000 Belgrade, V. Stepe 51., Serbia}}
\email{bsaric@ptt.yu}
\date{April 26, 2000}
\maketitle

\begin{abstract}
In the introduction part of this paper, first of all, the concept of
absolute integral sum of complex function is defined, as more general one
with respect to the concept of integral as well as of integral sum of
''ordinary'' integral calculus. In the main part of the paper, on the basis
of redefined concepts of a spatial derivative as well as of a residue of a
complex function, the fundamental results of \textit{Cauchy's} calculus of
residues of analytic functions are generalized and the concept of total
value of improper integral of complex function is defined. Finally, the
results in the form of fundamental lemmas, which are more general with
respect to the results of \textit{Jordan's} fundamental lemmas of \textit{%
Cauchy's} calculus of residues of analytic functions, are derived for the
general class of scalar complex function whose defining domain is a complex
plane.
\end{abstract}

\section{Introduction}

In the theory of complex analysis, more exactly of \textit{Cauchy's}
calculus of residues, the result formulated in the form of the following
theorem - \textit{Theorem 2, Subsection 3.1.1, Section 3.1, Chapter 3, p.
44, }\cite{M-K}

\begin{theorem}
\textit{If an analytic function }$f$\textit{\ has, in the extended complex
plane, only isolated singularities}\footnote{{\footnotesize The number of
isolated singularities of function }$f${\footnotesize \ must be finite,
because in the opposite case there exists the point of accumulation, which
is not isolated singularity.}}, \textit{then the sum of all its residues is
equal to zero}.$\blacktriangledown $
\end{theorem}

is one of fundamental ones.

The general case, in which the function $f$ has infinitely but a count of
many singularities, is an applicable, \cite{Ca}\textit{;} \cite{Pi} (taken
over from \cite{M-K}) and \cite{B-K}, and so according to that, it is
author's idea, presented in this article, to derive on the basis of the one
general theory being more general with respect to the theory of \textit{%
Cauchy's} calculus of residues, the results whose an importance and an
application are more general in comparison whit the fundamental results of 
\textit{Cauchy's} calculus of residues.

The function $f$, which is regular one in extended complex plane with the
exception at infinitely but a count of many points, can be said that is a
pseudo-analytic function, in view of the fact that it dose not belong to the
functional space of neither analytic nor non-analytic functions, \cite{M-K}.
Since the concept of residue of complex function, like that as it is
established for analytic functions, loses the sense for non-analytic as well
as for pseudo-analytic functions, then it is necessary to redefine the
concept of residue of complex function, more exactly to generalize it. The
concept of residue was first generalized by \textit{Poor}, \cite{Po1} and 
\cite{Po2} (taken over from \cite{M-K}) - \textit{Definition 1 and 2,
Subsection 2.2.2, Section 2.2, Chapter 2, pp. 38-39, }\cite{M-K}, for the
functional space of non-analytic functions.

In order to cam to the general definition of residue of a non-analytic
function, slightly general with the respect to \textit{Poor's} definitions,
it is necessary to introduce into the analysis the concept of complex
three-dimensional vector space $\vec{r}$\textit{:} $\vec{r}=x\vec{e}_{1}+iy%
\vec{e}_{2}+\varkappa \vec{n}$, of definite \textit{Euclidean} metric\textit{%
:} $ds^{2}=d\vec{r}\cdot d\vec{r}^{*}=dx^{2}+dy^{2}+d\varkappa ^{2}$ ($i$
denotes imaginary unit, the vector $\vec{r}^{*}$\textit{:} $\vec{r}^{*}=x%
\vec{e}_{1}-iy\vec{e}_{2}+\varkappa \vec{n}$, is a conjugate vector of a
position vector $\vec{r}$, and the vector $\vec{n}$ is an unit normal vector
of two-dimensional vector space $\vec{\varrho}$\textit{:} $\vec{\varrho}=x%
\vec{e}_{1}+iy\vec{e}_{2}$). In that case, the vector space $\vec{\varrho}$%
\textit{:} $\vec{\varrho}=x\vec{e}_{1}+iy\vec{e}_{2}$, is complex
two-dimensional vector plane of definite \textit{Euclidean} metric\textit{:} 
$dq^{2}=d\vec{\varrho}\cdot d\vec{\varrho}^{*}=dx^{2}+dy^{2}$. The map%
\textit{:} $z^{*}=x-iy$ and $z=x+iy$, is one-to-one map of the complex
two-dimensional vector planes $\vec{\varrho}$\textit{:} $\vec{\varrho}=z^{*}%
\vec{w}_{1}+z\vec{w}_{2}$ and $\vec{\varrho}=x\vec{e}_{1}+iy\vec{e}_{2}$%
\textit{;} $\vec{w}_{l}\cdot \vec{w}_{j}=\delta _{lj}$ ($\delta _{lj}$ is 
\textit{Kronecker's} $\delta $-symbol, more exactly an unit $2\times 2$
matrix), with \textit{Jacobian} of transformation\textit{:} $J=\left| 
\begin{array}{ll}
\frac{\partial z^{*}}{\partial x} & \frac{\partial z^{*}}{\partial iy} \\ 
\frac{\partial z}{\partial x} & \frac{\partial z}{\partial iy}
\end{array}
\right| =\left| 
\begin{array}{ll}
1 & -1 \\ 
1 & 1
\end{array}
\right| =2$. If the vector $\vec{r}_{A}$ is a position vector of an
arbitrary point $A$ of the vector space $\vec{r}$, then the set $\vec{r}_{g}$
of the points $\vec{r}_{A}$\textit{:} $\vec{r}_{g}=\left\{ \vec{r}_{A}\text{%
\textit{: }}A\in G\right\} $, defines some domain of three-dimensional
vector space $\vec{r}$. If the domain $\vec{r}_{g}$, bounded by a contour
surface $\vec{r}_{\gamma }$ of the vector space $\vec{r}$, is subdivided by
the planes being parallel to the co-ordinate planes, into $k_{1}$ elemental
sub-domains $\vartriangle _{j_{1}}\vec{r}_{g}$ bounded by the elemental
contour surfaces $\vartriangle _{j_{1}}\vec{r}_{\gamma }$\textit{\ (}$%
j_{1}=2,...,k_{1}$), then every sub-domain $\vartriangle _{j_{1}}\vec{r}_{g}$
of the domain $\vec{r}_{g}$ can be further subdivided into new sub-domains $%
\vartriangle _{j_{1},j_{2}}\vec{r}_{g}$ ($j_{2}=2,...,k_{2}$)\footnote{%
{\footnotesize If an infinite process of subdivision of domain }$\vec{r}_{g}$%
{\footnotesize ,\ is coming to each point }$\vec{r}_{N}${\footnotesize \ of
the domain }$\vec{r}_{g}$\textit{:}{\footnotesize \ }$\underset{n\rightarrow
+\infty }{\lim }\vartriangle _{j_{1},...,j_{n}}\vec{r}_{g}=\vec{r}_{N}$%
{\footnotesize , the subdivision }$D_{n}\vec{r}_{g}${\footnotesize \ is said
to be evenly spaced.}}. Let\textit{:} $\vartriangle _{j_{1},...,j_{n}}\sigma 
$ and $\vartriangle _{j_{1},...,j_{n}}v$, be measure numbers of an area of
an elemental contour surface $\vartriangle _{j_{1},...,j_{n}}\vec{r}_{\gamma
}$ as well as of a volume of an elemental sub-domain $\vartriangle
_{j_{1},...,j_{n}}\vec{r}_{g}$, respectively. In that case $d_{\vec{r}_{N}}v$%
\textit{:} $d_{\vec{r}_{N}}v=\underset{n\rightarrow +\infty }{\lim }%
\vartriangle _{j_{1},...,j_{n}}v=0$ and $d_{\vec{r}_{N}}\vec{\sigma}$\textit{%
:} $d_{\vec{r}_{N}}\vec{\sigma}=\underset{n\rightarrow +\infty }{\lim }\vec{n%
}_{\vec{r}_{\Delta _{j_{1},...,j_{n}}\vec{r}_{\gamma }}}\vartriangle
_{j_{1},...,j_{n}}\sigma =\vec{0}$ (the vector $\vec{n}_{\vec{r}_{\Delta
_{j_{1},...,j_{n}}\vec{r}_{\gamma }}}$ is an unit normal vector of an
elemental contour surface $\vartriangle _{j_{1},...,j_{n}}\vec{r}_{\gamma }$
at the point $\vec{r}_{\Delta _{j_{1},...,j_{n}}\vec{r}_{\gamma }}$), are an
infinitesimal volume and surface element, at the point $\vec{r}_{N}$ of the
domain $\vec{r}_{g}$, respectively.

Since the vectors $\vartriangle _{\vec{r}_{\Delta _{j_{1},...,j_{n}}\vec{r}%
_{\gamma }}}\vec{\sigma}$\textit{:} $\vartriangle _{\vec{r}_{\Delta
_{j_{1},...,j_{n}}\vec{r}_{\gamma }}}\vec{\sigma}=\vec{n}_{\vec{r}_{\Delta
_{j_{1},...,j_{n}}\vec{r}_{\gamma }}}\vartriangle _{j_{1},...,j_{n}}\sigma $%
, at an arbitrary point $\vec{r}_{\Delta _{j_{1},...,j_{n}}\vec{r}_{\gamma
}} $ of a part of an elemental contour surface $\vartriangle
_{j_{1},...,j_{n}}\vec{r}_{\gamma }$, which separates two elemental
sub-domains, have opposite orientations, then for every evenly spaced
subdivision $D_{n}\vec{r}_{g}$\textit{:} 
\begin{equation*}
D_{n}\vec{r}_{g}=\left\{ \vartriangle _{j_{1},...,j_{n}}\vec{r}_{g}\text{%
\textit{:} }j_{l}=2,...,k_{l}\,\left( l=1,2,...,n\right) \right\} ,
\end{equation*}
of the domain $\vec{r}_{g}$, bounded by the contour surface $\vec{r}_{\gamma
}$ of the vector space $\vec{r}$%
\begin{equation*}
\underset{n\rightarrow +\infty }{\lim }\overset{k_{1}}{\underset{j_{1}=2}{%
\sum }}\overset{k_{2}}{\underset{j_{2}=2}{\sum }}...\,\overset{k_{n}}{%
\underset{j_{n}=2}{\sum }}\vartriangle _{\vec{r}_{\Delta _{j_{1},...,j_{n}}%
\vec{r}_{\gamma }}}\vec{\sigma}=\underset{\vec{r}_{N}\in \vec{r}_{\gamma }}{%
\sum }d_{\vec{r}_{N}}\vec{\sigma}.
\end{equation*}

The infinite sum of zero vectors $\underset{\vec{r}_{N}\in \vec{r}_{\gamma }%
}{\sum }d_{\vec{r}_{N}}\vec{\sigma}$\textit{:} $\underset{\vec{r}_{N}\in 
\vec{r}_{\gamma }}{\sum }d_{\vec{r}_{N}}\vec{\sigma}=\vec{0}\infty $, as an
indefinite expression, in this acute case, reduces to the vector $\vec{P}$,
whose an intensity is equal to the area of the contour surface $\vec{r}%
_{\gamma }$\textit{:} $\underset{\vec{r}_{N}\in \vec{r}_{\gamma }}{\sum }d_{%
\vec{r}_{N}}\vec{\sigma}=\vec{P}$.

For an arbitrary scalar valued function $f\left( \vec{r}\right) $, defined
and bounded on the analyzed domain $\vec{r}_{g}$ of three-dimensional vector
space $\vec{r}$, it holds 
\begin{equation*}
\underset{n\rightarrow +\infty }{\lim }\overset{k_{1}}{\underset{j_{1}=2}{%
\sum }}\overset{k_{2}}{\underset{j_{2}=2}{\sum }}...\,\overset{k_{n}}{%
\underset{j_{n}=2}{\sum }}f\left( \vec{r}_{\Delta _{j_{1},...,j_{n}}\vec{r}%
_{\gamma }}\right) \vartriangle _{\vec{r}_{\Delta _{j_{1},...,j_{n}}\vec{r}%
_{\gamma }}}\vec{\sigma}=\underset{\vec{r}_{N}\in \vec{r}_{\gamma }}{\sum }%
f\left( \vec{r}_{N}\right) d_{\vec{r}_{N}}\vec{\sigma},
\end{equation*}
where $f\left( \vec{r}_{\Delta _{j_{1},...,j_{n}}\vec{r}_{\gamma }}\right) $
are values of the function $f\left( \vec{r}\right) $ at arbitrary points of
parts of elemental contour surfaces $\vartriangle _{j_{1},...,j_{n}}\vec{r}%
_{\gamma }$ which separate two elemental sub-domains, as well as at
arbitrary points of elemental contour surfaces $\vartriangle
_{j_{1},...,j_{n}}\vec{r}_{\gamma }$ which on the contour surface $\vec{r}%
_{\gamma }$. If the function $f\left( \vec{r}\right) $ is \textit{Riemann}%
-integrable over the domain $\vec{r}_{g}$, then for each evenly spaced
subdivision $D_{n}\vec{r}_{g}$ of domain $\vec{r}_{g}$ and for any a choice
of points $\vec{r}_{\Delta _{j_{1},...,j_{n}}\vec{r}_{\gamma }}$, there
exists a finite and an unique limiting value of the integral sum 
\begin{equation}
\underset{n\rightarrow +\infty }{\lim }\overset{k_{1}}{\underset{j_{1}=2}{%
\sum }}\overset{k_{2}}{\underset{j_{2}=2}{\sum }}...\,\overset{k_{n}}{%
\underset{j_{n}=2}{\sum }}f\left( \vec{r}_{\Delta _{j_{1},...,j_{n}}\vec{r}%
_{\gamma }}\right) \vartriangle _{\vec{r}_{\Delta _{j_{1},...,j_{n}}\vec{r}%
_{\gamma }}}\vec{\sigma}=\underset{\vec{r}_{\gamma }}{\overset{%
\circlearrowleft }{\iint }}f\left( \vec{r}\right) d\vec{\sigma}.  \label{1}
\end{equation}

Symbol $\underset{\vec{r}_{\gamma }}{\overset{\circlearrowleft }{\iint }}$
denotes an integration over the closed contour surface $\vec{r}_{\gamma }$,
in this case in the positive mathematical direction.

Now then, for a function $f\left( \vec{r}\right) $, which is \textit{Riemann}%
-integrable over the domain $\vec{r}_{g}$, the infinite sum of zero vectors%
\textit{:} $\underset{\vec{r}_{N}\in \vec{r}_{\gamma }}{\sum }f\left( \vec{r}%
_{N}\right) d_{\vec{r}_{N}}\vec{\sigma}=\infty \vec{0}$, as an indefinite
expression, is just equal to the integral of the function $f\left( \vec{r}%
\right) $%
\begin{equation}
\underset{\vec{r}_{N}\in \vec{r}_{\gamma }}{\sum }f\left( \vec{r}_{N}\right)
d_{\vec{r}_{N}}\vec{\sigma}=\underset{\vec{r}_{\gamma }}{\overset{%
\circlearrowleft }{\iint }}f\left( \vec{r}\right) d\vec{\sigma}.  \label{2}
\end{equation}

\begin{definition}
\textit{Absolute integral sums of a scalar valued function }$f\left( \vec{r}%
\right) $\textit{\ in the domain }$\vec{r}_{g}$\textit{\ bounded by contour
surface }$\vec{r}_{\gamma }$\textit{\ of the vector space }$\vec{r}$\textit{%
, are by definition} 
\begin{equation}
\underset{\vec{r}_{N}\in \vec{r}_{\gamma }}{\sum }f\left( \vec{r}_{N}\right)
d_{\vec{r}_{N}}\vec{\sigma},  \label{3}
\end{equation}
\begin{equation}
\underset{\vec{r}_{N}\in \vec{r}_{g}}{\sum }f\left( \vec{r}_{N}\right) d_{%
\vec{r}_{N}}v.\blacktriangledown  \label{4}
\end{equation}
\end{definition}

\begin{definition}
\textit{Absolute integral sums of a vector valued function }$\vec{F}\left( 
\vec{r}\right) $\textit{\ in the domain }$\vec{r}_{g}$\textit{\ bounded by
contour surface }$\vec{r}_{\gamma }$\textit{\ of the vector space }$\vec{r}$%
\textit{, are by definition} 
\begin{equation}
\underset{\vec{r}_{N}\in \vec{r}_{\gamma }}{\sum }\vec{F}\left( \vec{r}%
\right) \cdot d_{\vec{r}_{N}}\vec{\sigma},\underset{\vec{r}_{N}\in \vec{r}%
_{\gamma }}{\sum }d_{\vec{r}_{N}}\vec{\sigma}\times \vec{F}\left( \vec{r}%
\right) ,  \label{5}
\end{equation}
\begin{equation}
\underset{\vec{r}_{N}\in \vec{r}_{g}}{\sum }\vec{F}\left( \vec{r}\right) d_{%
\vec{r}_{N}}v.\blacktriangledown  \label{6}
\end{equation}
\end{definition}

\section{The main results}

\subsection{Spatial differentiability of a complex function}

On the basis of an integral equality of definition of spatial derivative of
scalar valued function $f\left( \vec{r}\right) $, in the \textit{Jung}
sense, \cite{Ju} (taken over from \cite{M-T}) - \textit{Definition 2,
Section 12.3, Chapter 12, p. 291, }\cite{M-T}\textit{\ }- for an arbitrary
uniform scalar valued function $f\left( \vec{r}\right) $ of the complex
three-dimensional vector space $\vec{r}$\textit{:} $\vec{r}=\vec{\varrho}%
+\varkappa \vec{n}$, as an ambient space of the complex two-dimensional
vector space $\vec{\varrho}$

\begin{definition}
\label{spatial}\textit{If and only if for every evenly spaced subdivision }$%
D_{n}\vec{r}_{g}$\textit{\ of the domain }$\vec{r}_{g}$\textit{\ bounded by
the contour surface }$\vec{r}_{\gamma }$\textit{\ of the vector space }$\vec{%
r}$\textit{\ and for every elemental sub-domain }$\vartriangle
_{j_{1},...,j_{n}}\vec{r}_{g}$\textit{\ of the domain }$\vec{r}_{g}$\textit{%
, the sequence of reduced absolute integral sums }$\vec{A}_{\vartriangle
_{j_{1},...,j_{n}}\vec{r}_{\gamma }}$\textit{:} 
\begin{equation}
\vec{A}_{\vartriangle _{j_{1},...,j_{n}}\vec{r}_{\gamma }}=\frac{1}{%
\vartriangle _{j_{1},...,j_{n}}v}\underset{\vec{r}_{N}\in \vartriangle
_{j_{1},...,j_{n}}\vec{r}_{\gamma }}{\sum }f\left( \vec{r}_{N}\right) d_{%
\vec{r}_{N}}\vec{\sigma},  \label{7}
\end{equation}
\textit{converges} 
\begin{equation}
\underset{\vec{r}_{\gamma }\rightarrow \vec{r}_{N}}{\lim }\frac{1}{V}%
\underset{\vec{r}_{N}\in \vec{r}_{\gamma }}{\sum }f\left( \vec{r}_{N}\right)
d_{\vec{r}_{N}}\vec{\sigma}=\underset{n\rightarrow +\infty }{\lim }\vec{A}%
_{\vartriangle _{j_{1},...,j_{n}}\vec{r}_{\gamma }}=\vec{A}_{\vec{r}_{N}},
\label{8}
\end{equation}
\textit{where }$V=\underset{\vec{r}_{N}\in \vec{r}_{g}}{\sum }d_{\vec{r}%
_{N}}v$\textit{, the function }$f\left( \vec{r}\right) $\textit{\ is a
spatial differentiable over the domain }$\vec{r}_{g}$\textit{.}$%
\blacktriangledown $
\end{definition}

The domain $\vec{r}_{g}$ of the vector space $\vec{r}$, such that at all
points of $\vec{r}_{g}$ the scalar valued function $f\left( \vec{r}\right) $
is a spatial differentiable, is a regular domain of the function $f\left( 
\vec{r}\right) .$ The points $\vec{r}_{N}$ of the vector space $\vec{r}$, at
which the function is not differentiable, are singular points of the
function $f\left( \vec{r}\right) $, and the domain $\vec{r}_{g}$, such that
the function $f\left( \vec{r}\right) $ is differentiable almost everywhere%
\footnote{{\footnotesize If the function }$f${\footnotesize \ possesses
someone feature everywhere with the exception, at most, of a set of points of%
} {\footnotesize \textit{Lebesgue's }measure zero, the function }$f$%
{\footnotesize \ is said to possess that feature almost everywhere.}} over $%
\vec{r}_{g}$, is a singular domain of the function $f\left( \vec{r}\right) $%
. The singular points $\vec{r}_{N}$ of the domain $\vec{r}_{g}$, at which
the function $f\left( \vec{r}\right) $ is bounded, are apparent singular
points of the function $f\left( \vec{r}\right) $. Singular domain $\vec{r}%
_{g}$, such that the function $f\left( \vec{r}\right) $ is bounded on $\vec{r%
}_{g}$, is an apparent singular domain of the function $f\left( \vec{r}%
\right) $.

If the function $f\left( \vec{r}\right) $ is spatial differentiable over the
domain $\vec{r}_{g}$, then on the basis of equality (\ref{7}) of \textit{%
Definition \ref{spatial}}, for an arbitrary evenly spaced subdivision $D_{n}%
\vec{r}_{g}$\textit{:} 
\begin{equation*}
D_{n}\vec{r}_{g}=\left\{ \vartriangle _{j_{1},...,j_{n}}\vec{r}_{g}\text{%
\textit{:} }j_{l}=2,...,k_{l}\,\left( l=1,2,...,n\right) \right\}
\end{equation*}
and for an arbitrary elemental sub-domain $\vartriangle _{j_{1},...,j_{n}}%
\vec{r}_{g}$ of the regular domain $\vec{r}_{g}$, it holds

\begin{equation}
\underset{\vec{r}_{N}\in \vartriangle _{j_{1},...,j_{n}}\vec{r}_{\gamma }}{%
\sum }f\left( \vec{r}_{N}\right) d_{\vec{r}_{N}}\vec{\sigma}=\vec{A}%
_{\vartriangle _{j_{1},...,j_{n}}\vec{r}_{\gamma }}\vartriangle
_{j_{1},...,j_{n}}v.  \label{9}
\end{equation}

On the one hand, having in view the fact that surface elements $d_{\vec{r}%
_{N}}\vec{\sigma}$ at each point $\vec{r}_{N}$ of the part of an elemental
contour surface $\vartriangle _{j_{1},...,j_{n}}\vec{r}_{\gamma }$, which
separates two elemental sub-domain are oppositely directed, for every level
of evenly spaced subdivision $D_{n}\vec{r}_{g}$ of the domain $\vec{r}_{g}$
it is obtained that 
\begin{equation*}
\overset{k_{1}}{\underset{j_{1}=2}{\sum }}\overset{k_{2}}{\underset{j_{2}=2}{%
\sum }}...\,\overset{k_{n}}{\underset{j_{n}=2}{\sum }}\underset{\vec{r}%
_{N}\in \vartriangle _{j_{1},...,j_{n}}\vec{r}_{\gamma }}{\sum }f\left( \vec{%
r}_{N}\right) d_{\vec{r}_{N}}\vec{\sigma}=\underset{\vec{r}_{N}\in \vec{r}%
_{\gamma }}{\sum }f\left( \vec{r}_{N}\right) d_{\vec{r}_{N}}\vec{\sigma}.
\end{equation*}

On the other hand, on the basis of a convergence of reduced absolute
integral sums, more exactly, of equality (\ref{8}) of \textit{Definition \ref
{spatial}}, it follows that 
\begin{equation*}
\underset{n\rightarrow +\infty }{\lim }\overset{k_{1}}{\underset{j_{1}=2}{%
\sum }}\overset{k_{2}}{\underset{j_{2}=2}{\sum }}...\,\overset{k_{n}}{%
\underset{j_{n}=2}{\sum }}\vec{A}_{\vartriangle _{j_{1},...,j_{n}}\vec{r}%
_{\gamma }}\vartriangle _{j_{1},...,j_{n}}v=\underset{\vec{r}_{N}\in \vec{r}%
_{g}}{\sum }\vec{A}_{\vec{r}_{N}}d_{\vec{r}_{N}}v.
\end{equation*}

Accordingly, it is finally obtained that 
\begin{equation}
\underset{\vec{r}_{N}\in \vec{r}_{\gamma }}{\sum }f\left( \vec{r}_{N}\right)
d_{\vec{r}_{N}}\vec{\sigma}=\underset{\vec{r}_{N}\in \vec{r}_{g}}{\sum }\vec{%
A}_{\vec{r}_{N}}d_{\vec{r}_{N}}v.  \label{10}
\end{equation}

If the function $f\left( \vec{r}\right) $ is defined on and is continuous on
the domain $\vec{r}_{g}$ bounded by the contour surface $\vec{r}_{\gamma }$
of the vector space $\vec{r}$, more exactly, is integrable over the domain $%
\vec{r}_{g}$, which is a regular domain of the function in the sense of 
\textit{Definition \ref{spatial}}, then on the basis of derived result (\ref
{10}) and of equality (\ref{2}), on the one hand 
\begin{equation}
\overset{\circlearrowleft }{\underset{\vec{r}_{\gamma }}{\iint }}f\left( 
\vec{r}\right) d\vec{\sigma}=\underset{\vec{r}_{N}\in \vec{r}_{g}}{\sum }%
\vec{A}_{\vec{r}_{N}}d_{\vec{r}_{N}}v,  \label{11}
\end{equation}
and on the other 
\begin{equation}
\underset{\vec{r}_{\gamma }\rightarrow \vec{r}_{N}}{\lim }\frac{1}{V}%
\overset{\circlearrowleft }{\underset{\vec{r}_{\gamma }}{\iint }}f\left( 
\vec{r}\right) d\vec{\sigma}=\vec{A}_{\vec{r}_{N}},  \label{12}
\end{equation}
where $V=\frac{1}{3}\overset{\circlearrowleft }{\underset{\vec{r}_{\gamma }}{%
\iint }}\vec{r}\cdot d\vec{\sigma}$.

If the vector valued function $\vec{A}\left( \vec{r}\right) $ ($\vec{A}%
\left( \vec{r}_{N}\right) =\vec{A}_{\vec{r}_{N}}$) is also defined on and is
continuous on the domain $\vec{r}_{g}$, more exactly, is integrable over the
domain $\vec{r}_{g}$, the equality\thinspace (\ref{11}) reduces to an
integral equality 
\begin{equation}
\overset{\circlearrowleft }{\underset{\vec{r}_{\gamma }}{\iint }}f\left( 
\vec{r}\right) d\vec{\sigma}=\underset{\vec{r}_{g}}{\iiint }\vec{A}\left( 
\vec{r}\right) dv.  \label{13}
\end{equation}

Clearly, the vector valued function $\vec{A}\left( \vec{r}\right) $ is a
function of spatial derivative of the continuous function $f\left( \vec{r}%
\right) $\textit{:} $\vec{A}\left( \vec{r}\right) =\nabla \cdot f\left( \vec{%
r}\right) $, and $\nabla $ is a \textit{Hamiltonian }operator of spatial
differentiability, \cite{M-T}.

In view of the fact that\textit{:} $d\vec{\sigma}=d\vec{\varrho}\times
d\varkappa \vec{n}+dz^{*}dz\vec{n}$ and $dv=d\vec{r}\cdot d\vec{\sigma}%
=dz^{*}dzd\varkappa $, for the function $f\left( \vec{r}\right) $ defined on
and continuous on (integrable over) the bounded domain $\vec{r}_{g}$ of the
vector space\thinspace $\vec{r}$, and which has defined and continuous
(integrable) partial derivatives on $\vec{r}_{g}$, from the integral
equality (\ref{13}) it follows that 
\begin{equation}
\overset{\circlearrowleft }{\underset{\vec{r}_{\gamma }}{\iint }}f\left( 
\vec{r}\right) dzd\varkappa =\underset{\vec{r}_{g}}{\iiint }\frac{\partial }{%
\partial z^{*}}f\left( \vec{r}\right) dz^{*}dzd\varkappa ,  \label{14}
\end{equation}
\begin{equation}
-\overset{\circlearrowleft }{\underset{\vec{r}_{\gamma }}{\iint }}f\left( 
\vec{r}\right) dz^{*}d\varkappa =\underset{\vec{r}_{g}}{\iiint }\frac{%
\partial }{\partial z}f\left( \vec{r}\right) dz^{*}dzd\varkappa ,  \label{15}
\end{equation}
\begin{equation}
\overset{\circlearrowleft }{\underset{\vec{r}_{\gamma }}{\iint }}f\left( 
\vec{r}\right) dz^{*}dz=\underset{\vec{r}_{g}}{\iiint }\frac{\partial }{%
\partial \varkappa }f\left( \vec{r}\right) dz^{*}dzd\varkappa .  \label{16}
\end{equation}

For the complex vector valued function $\vec{F}\left( \vec{r}\right) $%
\textit{:} $\vec{F}\left( \vec{r}\right) =P\left( \vec{r}\right) \vec{w}%
_{1}+Q\left( \vec{r}\right) \vec{w}_{2}+R\left( \vec{r}\right) \vec{n}$,
whose the components are defined and continuous (integrable) functions
having defined and continuous (integrable) partial derivatives on the
bounded domain $\vec{r}_{g}$ of the vector space\thinspace $\vec{r}$, on the
basis of the preceding derived results 
\begin{equation}
\overset{\circlearrowleft }{\underset{\vec{r}_{\gamma }}{\iint }}\vec{F}%
\left( \vec{r}\right) \cdot d\vec{\sigma}=\underset{\vec{r}_{g}}{\iiint }%
\nabla \cdot \vec{F}\left( \vec{r}\right) dv,  \label{17}
\end{equation}
\begin{equation}
\overset{\circlearrowleft }{\underset{\vec{r}_{\gamma }}{\iint }}d\vec{\sigma%
}\times \vec{F}\left( \vec{r}\right) =\underset{\vec{r}_{g}}{\iiint }\nabla
\times \vec{F}\left( \vec{r}\right) dv.  \label{18}
\end{equation}

The integral equalities\textit{:} (\ref{17}) and (\ref{18}), are analogous
to integral equalities of the well-known \textit{Gauss-Ostrogradski's} 
\textit{Theorem} for a real vector valued function $\vec{F}\left( \vec{r}%
\right) $ of the real three-dimensional vector space $\vec{r}$, \cite{M-T}.

In view of the fact that the limiting value $\vec{A}_{\vec{r}_{N}}$ of
sequence of reduced absolute integral sums $\vec{A}_{\vartriangle
_{j_{1},...,j_{n}}\vec{r}_{\gamma }}$, of equality (\ref{18}) of \textit{%
Definition \ref{spatial}} 
\begin{equation*}
\underset{\vec{r}_{\gamma }\rightarrow \vec{r}_{N}}{\lim }\frac{1}{V}%
\underset{\vec{r}_{N}\in \vec{r}_{\gamma }}{\sum }f\left( \vec{r}_{N}\right)
d_{\vec{r}_{N}}\vec{\sigma}=\underset{n\rightarrow +\infty }{\lim }\vec{A}%
_{\vartriangle _{j_{1},...,j_{n}}\vec{r}_{\gamma }}=\vec{A}_{\vec{r}_{N}},
\end{equation*}
dose not depend on the form of a contour surface $\vec{r}_{\gamma }$
bounding the domain $\vec{r}_{g}$ of the vector space $\vec{r}$, it follows
that if $\vec{r}_{s}\left( \vec{r}_{N},d\delta \right) $ is an
infinitesimally small spherical surface centred at $\vec{r}_{N}$ and of
radius $d\delta $, then for each point $\vec{r}_{N}$ lying inside $\vec{r}%
_{g}$ ($\vec{r}_{N}\in int.\vec{r}_{g}$, where $int.\vec{r}_{g}$ is an
interior of the domain $\vec{r}_{g}$) 
\begin{equation}
\vec{A}_{\vec{r}_{N}}=\frac{1}{d_{\vec{r}_{N}}v}\overset{\circlearrowleft }{%
\underset{\vec{r}_{s}\left( \vec{r}_{N},d\delta \right) }{\iint }}f\left( 
\vec{r}\right) d\vec{\sigma}.  \label{19}
\end{equation}

For a point $\vec{r}_{N}$ on the boundary $\vec{r}_{\gamma }$ of the domain $%
\vec{r}_{g}$ ($\vec{r}_{N}\in \vec{r}_{\gamma }$) 
\begin{equation}
\vec{A}_{\vec{r}_{N}}=\frac{1}{d_{\vec{r}_{N}}v}\overset{\circlearrowleft }{%
\underset{int.\vec{r}_{s}\left( \vec{r}_{N},d\delta \right) }{\iint }}%
f\left( \vec{r}\right) d\vec{\sigma},  \label{20}
\end{equation}
where $int.\vec{r}_{s}\left( \vec{r}_{N},d\delta \right) $ is a part of an
infinitesimally small spherical surface $\vec{r}_{s}\left( \vec{r}%
_{N},d\delta \right) $ that lies inside $\vec{r}_{g}$.

\begin{description}
\item[Comment]  The spacial derivative of a real valued function $f\left(
x\right) $ of one variable $x$, which is defined on the segment $\left[
a,b\right] $ of the real axis $R^{1}$, is by definition 
\begin{equation*}
\underset{\left[ a,b\right] \rightarrow c}{\lim }\frac{f\left( b\right)
-f\left( a\right) }{b-a}=A_{c},
\end{equation*}
more exactly, on the one hand, on the basis of relation (\ref{19}), for
point $c$ lying inside $\left[ a,b\right] $%
\begin{equation*}
A_{c}=\frac{1}{d_{c}x}\left[ f\left( c^{+}\right) -f\left( c^{-}\right)
\right] =\underset{2\vartriangle x\rightarrow d_{c}x}{\lim }\frac{f\left(
c+\vartriangle x\right) -f\left( c-\vartriangle x\right) }{2\vartriangle x},
\end{equation*}
and on the other, on the basis of relation (\ref{20}), for boundary points
of $\left[ a,b\right] $\textit{:} $a$ and $b$%
\begin{equation*}
A_{a}=\frac{1}{d_{a}x}\left[ f\left( a^{+}\right) -f\left( a\right) \right] =%
\underset{\vartriangle x\rightarrow d_{a}x}{\lim }\frac{f\left(
a+\vartriangle x\right) -f\left( a\right) }{\vartriangle x},
\end{equation*}
\begin{equation*}
A_{b}=\frac{1}{d_{b}x}\left[ f\left( b\right) -f\left( b^{-}\right) \right] =%
\underset{\vartriangle x\rightarrow d_{b}x}{\lim }\frac{f\left( b\right)
-f\left( b-\vartriangle x\right) }{\vartriangle x}.
\end{equation*}

According to that, as well as to the relation of equality (\ref{11}) 
\begin{equation*}
f\left( b\right) -f\left( a\right) =\underset{c\in \left[ a,b\right] }{\sum }%
A_{c}d_{c}x.
\end{equation*}

If the function $A\left( x\right) $ ($A\left( c\right) =A_{c}$) is an
integrable function over the segment $\left[ a,b\right] $, then the function 
$f\left( x\right) $ is a continuous function on $\left[ a,b\right] $ and 
\begin{equation*}
f\left( b\right) -f\left( a\right) =\overset{b}{\underset{a}{\int }}A\left(
x\right) dx,
\end{equation*}
where $A\left( x\right) =\nabla f\left( x\right) =\frac{df\left( x\right) }{%
dx}$.$\blacktriangledown $
\end{description}

\subsection{Residue of a complex function}

Based on the functional equality (\ref{8}) of \textit{Definition \ref
{spatial}}, more exactly, on the relations of equality\textit{:} (\ref{19})
and (\ref{20}), if a domain $\vec{r}_{g}$ bounded by the contour
surface\thinspace $\vec{r}_{\gamma }$ of the vector space $\vec{r}$, is a
regular domain of the function $f\left( \vec{r}\right) $, then at each point 
$\vec{r}_{N}$ of the domain $\vec{r}_{g}$\textit{:} $\vec{A}_{\vec{r}_{N}}d_{%
\vec{r}_{N}}v=\vec{0}$.

On the other hand, if a domain $\vec{r}_{g}$ bounded by the contour
surface\thinspace $\vec{r}_{\gamma }$ of the vector space $\vec{r}$, is a
singular domain of the function $f\left( \vec{r}\right) $, the absolute
integral sum of the function $\vec{A}\left( \vec{r}\right) $\textit{:} $%
\underset{\vec{r}_{N}\in int.\vec{r}_{g}}{\sum }\vec{A}_{\vec{r}_{N}}d_{\vec{%
r}_{N}}v$, can be subdivided into two absolute integral sums 
\begin{equation*}
\underset{\vec{r}_{N}\in \vec{r}_{g}}{\sum }\vec{A}_{\vec{r}_{N}}d_{\vec{r}%
_{N}}v=\underset{\vec{r}_{N}\in v.p.\vec{r}_{g}}{\sum }\vec{A}_{\vec{r}%
_{N}}d_{\vec{r}_{N}}v+\underset{\vec{r}_{N}\in s.p.\vec{r}_{g}}{\sum }\vec{A}%
_{\vec{r}_{N}}d_{\vec{r}_{N}}v,
\end{equation*}
where\textit{:} $v.p.\vec{r}_{g}$ and $s.p.\vec{r}_{g}$, are sets of the
regular and of the singular points of the function $f\left( \vec{r}\right) $
in the singular domain $\vec{r}_{g}$, respectively. Clearly, at all regular
points of the singular domain of the function $f\left( \vec{r}\right) $%
\textit{:} $\vec{A}_{\vec{r}_{N}}d_{\vec{r}_{N}}v=\vec{0}$, in other words
the first absolute integral sum of the function $\vec{A}\left( \vec{r}%
\right) $, on the left hand side of the preceding relation, is an infinite
sum of zero vectors. At each singular point $\vec{r}_{N}$ of the domain $%
\vec{r}_{g}$, at which the sequence of the reduced absolute integral sum $%
\vec{A}_{\vartriangle _{j_{1},...,j_{n}}\vec{r}_{\gamma }}$ definitely
diverges, $\vec{A}_{\vec{r}_{N}}d_{\vec{r}_{N}}v$ reduces to the indefinite
expression, more exactly, to either definite or indefinite vector value of
the extended vector space $\vec{r}\cup \vec{r}_{\infty }$, where $\vec{r}%
_{\infty }$ is a set of infinite points.

\subsubsection{The potential of a point with respect to a contour surface of
integration}

In the real two-dimensional vector space $\vec{\varrho}$\textit{:} $\vec{%
\varrho}=x\vec{e}_{1}+y\vec{e}_{2}$, an intensity of the vector $\frac{%
\left( \vec{\varrho}\times d\vec{\varrho}\right) }{\vec{\varrho}\cdot \vec{%
\varrho}}$\textit{:} $\frac{\left( \vec{\varrho}\times d\vec{\varrho}\right)
\cdot \vec{n}}{\vec{\varrho}\cdot \vec{\varrho}}=\frac{xdy-ydx}{x^{2}+y^{2}}%
=d\arctan \frac{y}{x}$, defines an intensity of an infinitesimal change of
the unit vector $\vec{\varrho}_{o}$ of a position vector $\vec{\varrho}$ of
an arbitrary point of two-dimensional vector space $\vec{\varrho}$ with
respect to the origin\textit{:} 
\begin{equation*}
\frac{\left( \vec{\varrho}\times d\vec{\varrho}\right) \cdot \vec{n}}{\vec{%
\varrho}\cdot \vec{\varrho}}=\left| d\vec{\varrho}_{o}\right| =d\varphi .
\end{equation*}

In the case of complex two-dimensional vector space $\vec{\varrho}$\textit{:}
$\vec{\varrho}=\left| \vec{\varrho}\right| \vec{\varrho}_{o}$\footnote{%
{\footnotesize On the one hand} 
\begin{equation*}
e^{2i\arctan \frac{y}{x}}=\cos \left( 2\arctan \frac{y}{x}\right) +i\sin
\left( 2\arctan \frac{y}{x}\right) =
\end{equation*}
\begin{equation*}
=\left[ 1+i\tan \left( 2\arctan \frac{y}{x}\right) \right] \cos \left(
2\arctan \frac{y}{x}\right) =
\end{equation*}
\begin{equation*}
=\left[ 1+i\tan \left( \arctan \frac{y}{x}\right) \right] ^{2}\cos
^{2}\left( \arctan \frac{y}{x}\right) =
\end{equation*}
\begin{equation*}
=\left( 1+i\frac{y}{x}\right) ^{2}\cos ^{2}\left( \arctan \frac{y}{x}\right)
,
\end{equation*}
{\footnotesize and on the other} 
\begin{equation*}
\cos ^{2}\left( \arctan \frac{y}{x}\right) =1-\sin ^{2}\left( \arctan \frac{y%
}{x}\right) ,
\end{equation*}
{\footnotesize more exactly} 
\begin{equation*}
\cos ^{2}\left( \arctan \frac{y}{x}\right) \left[ 1+\tan ^{2}\left( \arctan 
\frac{y}{x}\right) \right] =\cos ^{2}\left( \arctan \frac{y}{x}\right)
\left[ 1+\left( \frac{y}{x}\right) ^{2}\right] =1{\footnotesize .}
\end{equation*}
\par
{\footnotesize From preceding equalities it follows that }$e^{2i\arctan 
\frac{y}{x}}=\frac{\left( 1+i\frac{y}{x}\right) ^{2}}{1+\left( \frac{y}{x}%
\right) ^{2}}=\frac{x+iy}{x-iy}.$ {\footnotesize Accordingly,} 
{\footnotesize since for }$\varphi =\arctan \frac{y}{x}$\textit{:}%
{\footnotesize \ }$z=z^{*}e^{2i\varphi }${\footnotesize , then }$\vec{\varrho%
}=\sqrt{\frac{\vec{\varrho}\cdot \vec{\varrho}^{*}}{2}}\left( e^{-i\varphi }%
\vec{w}_{1}+e^{i\varphi }\vec{w}_{2}\right) =\left| \vec{\varrho}\right| 
\vec{\varrho}_{o}${\footnotesize , where\textit{:} }$\left| \vec{\varrho}%
\right| =\sqrt{\vec{\varrho}\cdot \vec{\varrho}^{*}}${\footnotesize \ and }$%
\vec{\varrho}_{o}=\frac{1}{\sqrt{2}}\left( e^{-i\varphi }\vec{w}%
_{1}+e^{i\varphi }\vec{w}_{2}\right) ${\footnotesize .}}, it follows that 
\begin{equation*}
\overset{\circlearrowleft }{\underset{\vec{\varrho}_{\gamma }}{\int }}\frac{%
\left( \vec{\varrho}\times d\vec{\varrho}\right) \cdot \vec{n}}{\vec{\varrho}%
\cdot \vec{\varrho}^{*}}=\overset{\circlearrowleft }{\underset{\vec{\varrho}%
_{\gamma }}{\int }}\left( \vec{\varrho}_{o}\times d\vec{\varrho}_{o}\right)
\cdot \vec{n}=i\overset{\circlearrowleft }{\underset{\vec{\varrho}_{\gamma }%
}{\int }}d\varphi .
\end{equation*}

\begin{definition}
\textit{A potential }$p_{\vec{\varrho}_{\gamma }\leftrightarrow \vec{\varrho}%
_{N}}$\textit{\ of a point }$\vec{\varrho}_{N}$\textit{\ with respect to a
contour }$\vec{\varrho}_{\gamma }$\textit{\ bounding the domain }$\vec{%
\varrho}_{g}$\textit{\ of the complex plane }$\vec{\varrho}$\textit{, is by
definition} 
\begin{equation}
p_{\vec{\varrho}_{\gamma }\leftrightarrow \vec{\varrho}_{N}}=\overset{%
\circlearrowleft }{\underset{\vec{\varrho}_{\gamma }}{\int }}\frac{\left[
\left( \vec{\varrho}-\vec{\varrho}_{N}\right) \times d\vec{\varrho}\right]
\cdot \vec{n}}{\left( \vec{\varrho}-\vec{\varrho}_{N}\right) \cdot \left( 
\vec{\varrho}-\vec{\varrho}_{N}\right) ^{*}}.\blacktriangledown  \label{21}
\end{equation}
\end{definition}

In view of the fact that differential $d\varphi $ is an absolute one, then
for every closed path of an integration $\vec{\varrho}_{\gamma }$, by which
a domain $\vec{\varrho}_{g}$ of the complex plane $\vec{\varrho}$ is
bounded, it holds 
\begin{equation}
p_{\vec{\varrho}_{\gamma }\leftrightarrow \vec{\varrho}_{N}}=i\overset{%
\circlearrowleft }{\underset{\vec{\varrho}_{\gamma }}{\int }}d\theta
=\left\{ 
\begin{array}{l}
2\pi i,\vec{\varrho}_{N}\in int.\,\vec{\varrho}_{g} \\ 
0,\vec{\varrho}_{N}\notin \vec{\varrho}_{g}
\end{array}
\right\} ,  \label{22}
\end{equation}
where $int.\,\vec{\varrho}_{g}$ is an interior of the domain $\vec{\varrho}%
_{g}$.

In the case in which a point $\vec{\varrho}_{N}$ belongs to the boundary $%
\vec{\varrho}_{\gamma }$ of the domain $\vec{\varrho}_{g}$, the potential $%
p_{\vec{\varrho}_{\gamma }\leftrightarrow \vec{\varrho}_{N}}$ of the point $%
\vec{\varrho}_{N}$ with respect to $\vec{\varrho}_{\gamma }$ is defined to
be the sum of limiting values of an integral $i\overset{\circlearrowleft }{%
\underset{\vec{\varrho}_{\gamma }}{\int }}d\theta $ over the part of the
path of an integration $\vec{\varrho}_{\gamma }$ from the point $\vec{\varrho%
}_{A}$ to the point $\vec{\varrho}_{B}$ ($\vec{\varrho}_{A}$ and $\vec{%
\varrho}_{B}$ are intersection points of the path of an integration $\vec{%
\varrho}_{\gamma }$ and of some arbitrary small circle $\vec{\varrho}%
_{\delta }$ centred at $\vec{\varrho}_{N}$ and of radius $\delta $) and over
the circular arcs from the point $\vec{\varrho}_{B}$ to the point $\vec{%
\varrho}_{A}$, as the radius $\delta $ of the circle $\vec{\varrho}_{\delta
} $ tends to zero, in other words as the boundary points of the circular arcs%
\textit{:} $\vec{\varrho}_{A}$ and $\vec{\varrho}_{B}$, along the path of an
integration $\vec{\varrho}_{\gamma }$, tends to the point $\vec{\varrho}_{N}$%
.

Considering the fact that the limiting value of an integral $i\overset{%
\circlearrowleft }{\underset{\vec{\varrho}_{\gamma }}{\int }}d\theta $ over
the part of the path of an integration $\vec{\varrho}_{\gamma }$ from the
point $\vec{\varrho}_{A}$ to the point $\vec{\varrho}_{B}$ is equal to%
\textit{: }$\underset{\delta \rightarrow 0^{+}}{\lim }i\underset{\vec{\varrho%
}_{\gamma }}{\overset{\overset{\curvearrowleft }{\vec{\varrho}_{B}\vec{%
\varrho}_{A}}}{\int }}d\theta =$ $i\alpha $, and over the circular arcs%
\textit{:} $\underset{\delta \rightarrow 0^{+}}{\lim }i\underset{int.\vec{%
\varrho}_{\delta }}{\overset{\overset{\curvearrowright }{\vec{\varrho}_{B}%
\vec{\varrho}_{A}}}{\int }}d\theta =-i\alpha $ and $\underset{\delta
\rightarrow 0^{+}}{\lim }i\underset{ext.\vec{\varrho}_{\delta }}{\overset{%
\overset{\curvearrowleft }{\vec{\varrho}_{A}\vec{\varrho}_{B}}}{\int }}%
d\theta =i\left( 2\pi -\alpha \right) $, where\textit{:} $int.\vec{\varrho}%
_{\delta }$ and $ext.\vec{\varrho}_{\delta }$, are the circular arcs inside
and outside the domain $\vec{\varrho}_{g}$ respectively, and the angle $%
\alpha $ is a limit angle of tangent lines to the path of an integration $%
\vec{\varrho}_{\gamma }$ at the points\textit{:} $\vec{\varrho}_{A}$ and $%
\vec{\varrho}_{B}$, as the boundary points of the circular arcs\textit{:} $%
\vec{\varrho}_{A}$ and $\vec{\varrho}_{B}$, along the path of an integration 
$\vec{\varrho}_{\gamma }$, tends to the point $\vec{\varrho}_{N}$, in this
emphasized case 
\begin{equation}
p_{\vec{\varrho}_{\gamma }\leftrightarrow \vec{\varrho}_{N}}=\left\{ 
\begin{array}{l}
0 \\ 
2\pi i
\end{array}
\right\} .  \label{23}
\end{equation}

\begin{definition}
\label{pot}\textit{A potential }$p_{\vec{r}_{\gamma }\leftrightarrow \vec{r}%
_{N}}$\textit{\ of a point }$\vec{r}_{N}$\textit{\ with respect to a contour
surface }$\vec{r}_{\gamma }$\textit{\ bounding the domain }$\vec{r}_{g}$%
\textit{\ of the vector space }$\vec{r}$\textit{:} $\vec{r}=\vec{\varrho}%
+\varkappa \vec{n}$\textit{, is by definition} 
\begin{equation}
p_{\vec{r}_{\gamma }\leftrightarrow \vec{r}_{N}}=2p_{\vec{\varrho}_{\gamma
}\leftrightarrow \vec{\varrho}_{N}},  \label{24}
\end{equation}
\textit{where }$\vec{\varrho}_{\gamma }$\textit{: }$\vec{\varrho}_{\gamma }=%
\vec{r}_{\gamma }\cap \vec{\varrho}$\textit{, and }$\vec{\varrho}$\textit{\
is any complex plane such that }$\vec{r}_{N}=\vec{\varrho}_{N}$\textit{.}$%
\blacktriangledown $
\end{definition}

If one takes into consideration the fact that except the point $\vec{r}_{N}$%
, which is an inner point with respect to an infinitesimally small spherical
surface $\vec{r}_{s}\left( \vec{r}_{N},d\delta \right) $, the all remaining
points of the vector space $\vec{r}$ are external points, then according to
the defined concept of the potential $p_{\vec{r}_{\gamma }\leftrightarrow 
\vec{r}_{N}}$ of a point $\vec{r}_{N}$ with respect to a contour surface $%
\vec{r}_{\gamma }$ bounding a certain domain\thinspace $\vec{r}_{g}$ of the
vector space $\vec{r}$ and on the basis of \textit{Poor's} definition of
non-analytic function residue - \textit{Definition 1, Subsection 2.2.2,
Section 2.2, Chapter 2, p. 38,} \cite{M-K}

\begin{definition}
\label{res}\textit{A residue (}$Res$\textit{) of a scalar valued function }$%
f\left( \vec{r}\right) $\textit{\ at a point }$\vec{r}_{N}$\textit{\ of the
vector space\thinspace }$\vec{r}$\textit{, is by definition} 
\begin{equation}
\overset{\circlearrowleft }{\underset{\vec{r}_{s}\left( \vec{r}_{N},d\delta
\right) }{\iint }}f\left( \vec{r}\right) d\vec{\sigma}=p_{\vec{r}_{s}\left( 
\vec{r}_{N},d\delta \right) \leftrightarrow \vec{r}_{N}}\underset{\vec{r}=%
\vec{r}_{N}}{\overset{\rightarrow }{Res}}f\left( \vec{r}\right)
.\blacktriangledown  \label{25}
\end{equation}
\end{definition}

\begin{definition}
\label{res inf.}\textit{A residue (}$Res$\textit{) of a scalar valued
function }$f\left( \vec{r}\right) $\textit{\ at the set of the infinite
points }$\vec{r}_{\infty }$\textit{, is by definition} 
\begin{equation}
\underset{\vec{r}=\vec{r}_{\infty }}{\overset{\rightarrow }{Res}}f\left( 
\vec{r}\right) =-\underset{\vec{r}_{N}\in \vec{r}}{\sum }\underset{\vec{r}=%
\vec{r}_{N}}{\overset{\rightarrow }{Res}}f\left( \vec{r}\right)
.\blacktriangledown  \label{26}
\end{equation}
\end{definition}

On the one hand, based on the equality (\ref{25}) of \textit{Definition \ref
{res}}, at all points inside a certain singular domain $\vec{r}_{g}$ bounded
by the contour surface $\vec{r}_{\gamma }$ of the vector space $\vec{r}$ ($%
\vec{r}_{N}\in int.\vec{r}_{g}$) 
\begin{equation}
\vec{A}_{\vec{r}_{N}}d_{\vec{r}_{N}}v=\nabla _{\vec{r}_{N}}f\left( \vec{r}%
\right) d_{\vec{r}_{N}}v=p_{\vec{r}_{\gamma }\leftrightarrow \vec{r}_{N}}%
\underset{\vec{r}=\vec{r}_{N}}{\overset{\rightarrow }{Res}}f\left( \vec{r}%
\right) ,  \label{27}
\end{equation}
and on the other, at the points onto the boundary\thinspace $\vec{r}_{\gamma
}$ of the domain $\vec{r}_{g}$ ($\vec{r}_{N}\in \vec{r}_{\gamma }$) 
\begin{equation}
\vec{A}_{\vec{r}_{N}}d_{\vec{r}_{N}}v=p_{int.\vec{r}_{s}\left( \vec{r}%
_{N},d\delta \right) \leftrightarrow \vec{r}_{N}}\underset{\vec{r}=\vec{r}%
_{N}}{\overset{\rightarrow }{Res}}f\left( \vec{r}\right) ,  \label{28}
\end{equation}
where, according to equality (\ref{24})\textit{:} $p_{int.\vec{r}_{s}\left( 
\vec{r}_{N},d\delta \right) \leftrightarrow \vec{r}_{N}}=2p_{int.\vec{\varrho%
}_{s}\left( \vec{\varrho}_{N},d\delta \right) \leftrightarrow \vec{\varrho}%
_{N}}$ and 
\begin{equation*}
p_{int.\vec{\varrho}_{s}\left( \vec{\varrho}_{N},d\delta \right)
\leftrightarrow \vec{\varrho}_{N}}=i\underset{int.\vec{\varrho}_{s}\left( 
\vec{\varrho}_{N},d\delta \right) }{\overset{\curvearrowleft }{\int }}%
d\theta =\underset{\delta \rightarrow 0^{+}}{\lim }i\underset{int.\vec{%
\varrho}_{\delta }}{\overset{\overset{\curvearrowleft }{\vec{\varrho}_{B}%
\vec{\varrho}_{A}}}{\int }}d\theta =i\alpha ,
\end{equation*}
and $\alpha $ is an angle of tangent lines at the point $\vec{\varrho}_{N}$
onto the boundary $\vec{\varrho}_{\gamma }$ of the domain $\vec{\varrho}_{g}$%
\textit{:} $\vec{\varrho}_{g}=\vec{r}_{g}\cap \vec{\varrho}$ ($\vec{\varrho}%
_{\gamma }=\vec{r}_{\gamma }\cap \vec{\varrho}$).

Clearly 
\begin{equation*}
p_{ext.\vec{\varrho}_{s}\left( \vec{\varrho}_{N},d\delta \right)
\leftrightarrow \vec{\varrho}_{N}}=i\underset{ext.\vec{\varrho}_{s}\left( 
\vec{\varrho}_{N},d\delta \right) }{\overset{\curvearrowleft }{\int }}%
d\theta =\underset{\delta \rightarrow 0^{+}}{\lim }i\underset{ext.\vec{%
\varrho}_{\delta }}{\overset{\overset{\curvearrowleft }{\vec{\varrho}_{A}%
\vec{\varrho}_{B}}}{\int }}d\theta =i\left( 2\pi -\alpha \right) .
\end{equation*}

\begin{definition}
\label{v.p.}\textit{Cauchy's principal value (}$v.p.$\textit{) of an
improper integral of vector valued function }$\nabla f\left( \vec{r}\right) $%
\textit{\ on a certain domain }$\vec{r}_{g}$\textit{\ bounded by the contour
surface\thinspace }$\vec{r}_{\gamma }$\textit{\ of the vector space }$\vec{r}
$\textit{, is by definition} 
\begin{equation}
v.p.\underset{\vec{r}_{g}}{\iiint }\nabla f\left( \vec{r}\right) dv=%
\underset{\vec{r}_{N}\in v.p.\vec{r}_{g}}{\sum }\vec{A}_{\vec{r}_{N}}d_{\vec{%
r}_{N}}v,  \label{29}
\end{equation}
\textit{where }$v.p.\vec{r}_{g}$\textit{\ is a set of the regular points }$%
\vec{r}_{N}$\textit{\ of the function on the domain }$\vec{r}_{g}$\textit{.}$%
\blacktriangledown $
\end{definition}

\begin{definition}
\label{v.s.}\textit{Jordan's singular value (}$v.s.$\textit{) of an improper
integral of vector valued function }$\nabla f\left( \vec{r}\right) $\textit{%
\ on a certain domain }$\vec{r}_{g}$\textit{\ bounded by the contour
surface\thinspace }$\vec{r}_{\gamma }$\textit{\ of the vector space }$\vec{r}
$\textit{, is by definition} 
\begin{equation}
v.s.\underset{\vec{r}_{g}}{\iiint }\nabla f\left( \vec{r}\right) dv=%
\underset{\vec{r}_{N}\in v.s.\vec{r}_{g}}{\sum }\vec{A}_{\vec{r}_{N}}d_{\vec{%
r}_{N}}v,  \label{30}
\end{equation}
\textit{where }$v.p.\vec{r}_{g}$\textit{\ is a set of the singular points }$%
\vec{r}_{N}$\textit{\ of the function on the domain }$\vec{r}_{g}$\textit{.}$%
\blacktriangledown $
\end{definition}

\begin{definition}
\label{v.t.}\textit{The sum of Cauchy's principal value (}$v.p.$\textit{)
and of Jordan's singular value (}$v.s.$\textit{) of an improper integral is
a total value (}$v.t.$\textit{) of that improper integral.}$%
\blacktriangledown $
\end{definition}

Since the derived equality (\ref{10}) is also valid in the case in which the
sequences of reduced absolute integral sums $\vec{A}_{\vartriangle
_{j_{1},...,j_{n}}\vec{r}_{\gamma }}$ diverge, in other words it is also
valid for singular domains\thinspace $\vec{r}_{g}$ of the function $f\left( 
\vec{r}\right) $, then if the function $f\left( \vec{r}\right) $ is
integrable function over contour surface $\vec{r}_{\gamma }$ bounding a
certain singular domain\thinspace $\vec{r}_{g}$ of the function $f\left( 
\vec{r}\right) $, it follows that 
\begin{equation*}
\overset{\circlearrowleft }{\underset{\vec{r}_{\gamma }}{\iint }}f\left( 
\vec{r}\right) d\vec{\sigma}=\underset{\vec{r}_{N}\in v.p.\vec{r}_{g}}{\sum }%
\vec{A}_{\vec{r}_{N}}d_{\vec{r}_{N}}v+\underset{\vec{r}_{N}\in s.p.\vec{r}%
_{g}}{\sum }\vec{A}_{\vec{r}_{N}}d_{\vec{r}_{N}}v.
\end{equation*}

Based on the equalities\textit{:} (\ref{29}) and (\ref{30}), on the one hand
and on \textit{Definition \ref{v.t.}}, on the other hand, finally it is
obtained that 
\begin{equation}
\overset{\circlearrowleft }{\underset{\vec{r}_{\gamma }}{\iint }}f\left( 
\vec{r}\right) d\vec{\sigma}-v.p.\underset{\vec{r}_{g}}{\iiint }\nabla
f\left( \vec{r}\right) dv=\underset{\vec{r}_{N}\in v.s.\vec{r}_{g}}{\sum }p_{%
\vec{r}_{\gamma }\leftrightarrow \vec{r}_{N}}\underset{\vec{r}=\vec{r}_{N}}{%
\overset{\rightarrow }{Res}}f\left( \vec{r}\right) ,  \label{31}
\end{equation}
where $v.s.\vec{r}_{g}$ is a set of the singular points $\vec{r}_{N}$ of the
function $f\left( \vec{r}\right) $ on the singular domain $\vec{r}_{g}$
bounded by the contour surface $\vec{r}_{\gamma }$ of the vector space $\vec{%
r}$, more exactly 
\begin{equation}
\overset{\circlearrowleft }{\underset{\vec{r}_{\gamma }}{\iint }}f\left( 
\vec{r}\right) d\vec{\sigma}=v.t.\underset{\vec{r}_{g}}{\iiint }\nabla
f\left( \vec{r}\right) dv.  \label{32}
\end{equation}

If the function $f\left( \vec{r}\right) $ is not integrable over a contour
surface\thinspace $\vec{r}_{\gamma }$ bounding a certain singular domain $%
\vec{r}_{g}$ of the function, in other words if singularities of the
function $f\left( \vec{r}\right) $ lie not only inside but also onto contour
surface $\vec{r}_{\gamma }$, on the one hand 
\begin{equation*}
\underset{\vec{r}_{N}\in v.p.\vec{r}_{\gamma }}{\sum }f\left( \vec{r}%
_{N}\right) d_{\vec{r}_{N}}\vec{\sigma}+\underset{\vec{r}_{N}\in v.s.\vec{r}%
_{\gamma }}{\sum }\underset{int.\vec{r}_{s}\left( \vec{r}_{N},d\delta
\right) }{\overset{\curvearrowright }{\iint }}f\left( \vec{r}\right) d\vec{%
\sigma}=
\end{equation*}
\begin{equation*}
=v.p.\underset{\vec{r}_{g}}{\iiint }\nabla f\left( \vec{r}\right) dv+%
\underset{\vec{r}_{N}\in v.s.int.\vec{r}_{g}}{\sum }4\pi i\underset{\vec{r}=%
\vec{r}_{N}}{\overset{\rightarrow }{Res}}f\left( \vec{r}\right) ,
\end{equation*}
and on the other 
\begin{equation*}
\underset{\vec{r}_{N}\in v.p.\vec{r}_{\gamma }}{\sum }f\left( \vec{r}%
_{N}\right) d_{\vec{r}_{N}}\vec{\sigma}+\underset{\vec{r}_{N}\in v.s.\vec{r}%
_{\gamma }}{\sum }\underset{ext.\vec{r}_{s}\left( \vec{r}_{N},d\delta
\right) }{\overset{\curvearrowleft }{\iint }}f\left( \vec{r}\right) d\vec{%
\sigma}=
\end{equation*}
\begin{equation*}
=v.p.\underset{\vec{r}_{g}}{\iiint }\nabla f\left( \vec{r}\right) dv+%
\underset{\vec{r}_{N}\in v.s.\vec{r}_{g}}{\sum }4\pi i\underset{\vec{r}=\vec{%
r}_{N}}{\overset{\rightarrow }{Res}}f\left( \vec{r}\right) ,
\end{equation*}
more exactly 
\begin{equation}
v.t.\overset{\circlearrowleft }{\underset{\vec{r}_{\gamma }}{\iint }}f\left( 
\vec{r}\right) d\vec{\sigma}-v.p.\underset{\vec{r}_{g}}{\iiint }\nabla
f\left( \vec{r}\right) dv=\underset{\vec{r}_{N}\in v.s.\vec{r}_{g}}{\sum }p_{%
\vec{r}_{\gamma }\leftrightarrow \vec{r}_{N}}\underset{\vec{r}=\vec{r}_{N}}{%
\overset{\rightarrow }{Res}}f\left( \vec{r}\right) ,  \label{33}
\end{equation}
where 
\begin{equation*}
v.t.\overset{\circlearrowleft }{\underset{\vec{r}_{\gamma }}{\iint }}f\left( 
\vec{r}\right) d\vec{\sigma}=v.p.\overset{\circlearrowleft }{\underset{\vec{r%
}_{\gamma }}{\iint }}f\left( \vec{r}\right) d\vec{\sigma}+v.s.\overset{%
\circlearrowleft }{\underset{\vec{r}_{\gamma }}{\iint }}f\left( \vec{r}%
\right) d\vec{\sigma}=
\end{equation*}
\begin{equation*}
=\underset{\vec{r}_{N}\in v.p.\vec{r}_{\gamma }}{\sum }f\left( \vec{r}%
_{N}\right) d_{\vec{r}_{N}}\vec{\sigma}+\underset{\vec{r}_{N}\in v.s.\vec{r}%
_{\gamma }}{\sum }\left\{ 
\begin{array}{l}
-p_{int.\vec{r}_{s}\left( \vec{r}_{N},d\delta \right) \leftrightarrow \vec{r}%
_{N}} \\ 
p_{ext.\vec{r}_{s}\left( \vec{r}_{N},d\delta \right) \leftrightarrow \vec{r}%
_{N}}
\end{array}
\right\} \underset{\vec{r}=\vec{r}_{N}}{\overset{\rightarrow }{Res}}f\left( 
\vec{r}\right) .
\end{equation*}

\begin{description}
\item[Comment]  For a real valued function $f\left( x\right) $ of one
variable $x$, which is spatial differentiable almost everywhere over the
segment $\left[ a,b\right] $ of the real axis $R^{1}$ and defined at
boundary points\textit{:} $a$ and $b$, of the segment $\left[ a,b\right] $,
and on the basis of the result (\ref{31}) 
\begin{equation*}
f\left( b\right) -f\left( a\right) -v.p.\overset{b}{\underset{a}{\int }}%
\nabla f\left( x\right) dx=\underset{c\in v.s.\left[ a,b\right] }{\sum }%
p_{a,b\leftrightarrow c}\underset{x=c}{Res}f\left( x\right) ,
\end{equation*}
where\textit{:} $2p_{a,b\leftrightarrow c}=p_{\vec{\varrho}_{\gamma
}\leftrightarrow \vec{\varrho}_{c}}$ and $\vec{\varrho}_{\gamma }$\textit{:} 
$\vec{\varrho}_{\gamma }\cap R^{1}=\left\{ a,b\right\} $, as well as 
\begin{equation*}
p_{a,b\leftrightarrow c}\underset{x=c}{Res}f\left( x\right) =f\left(
c^{+}\right) -f\left( c^{-}\right) =A_{c}d_{c}x,
\end{equation*}
more exactly 
\begin{equation*}
p_{a,b\leftrightarrow c}\underset{x=c}{Res}f\left( x\right) =\underset{%
2\vartriangle x\rightarrow d_{c}x}{\lim }\left[ f\left( c+\vartriangle
x\right) -f\left( c-\vartriangle x\right) \right] .
\end{equation*}

Based on \textit{Definitions \ref{v.p.}, \ref{v.s.}} and \textit{\ref{v.t.}} 
\begin{equation*}
f\left( b\right) -f\left( a\right) =v.t.\overset{b}{\underset{a}{\int }}%
\nabla f\left( x\right) dx.\blacktriangledown
\end{equation*}
\end{description}

\begin{example}
The scalar valued function $f\left( x\right) =\log x$, where $\log $ denotes
principal logarithm, is spatial differentiable at all points of the segment $%
\left[ -a,b\right] $ of the real axis\thinspace $R^{1}$ ($a,b\in R_{+}^{1}$)
except at the point $x=0$. Since, on the one hand 
\begin{equation*}
p_{-a,b\leftrightarrow 0}\underset{x=0}{Res}f\left( x\right) =\underset{%
2\vartriangle x\rightarrow d_{0}x}{\lim }\left[ \log \left( \vartriangle
x\right) -\log \left( -\vartriangle x\right) \right] =\mp \pi i,
\end{equation*}
and on the other 
\begin{equation*}
v.p.\overset{b}{\underset{-a}{\int }}\nabla f\left( x\right) dx=\log \frac{b%
}{a},
\end{equation*}
it follows that $\log \frac{b}{a}\mp \pi i=v.t.\overset{b}{\underset{-a}{%
\int }}\frac{dx}{x}.\blacktriangledown $
\end{example}

\begin{example}
The scalar valued function $f\left( x\right) =x^{-1}$ is spatial
differentiable at all points of the segment $\left[ -a,b\right] $ of the
real axis\thinspace $R^{1}$ ($a,b\in R_{+}^{1}$) except at the point $x=0$.
Since, on the one hand 
\begin{equation*}
p_{-a,b\leftrightarrow 0}\underset{x=0}{Res}f\left( x\right) =\underset{%
2\vartriangle x\rightarrow d_{0}x}{\lim }\left[ \frac{1}{\vartriangle x}+%
\frac{1}{\vartriangle x}\right] =+\infty ,
\end{equation*}
and on the other 
\begin{equation*}
v.p.\overset{b}{\underset{-a}{\int }}\nabla f\left( x\right) dx=-\infty ,
\end{equation*}
then in this case the total value of an improper integral\textit{:} $%
\overset{b}{\underset{-a}{\int }}x^{-2}dx$, as an indefinite expression of
difference of infinities, has exactly definite value 
\begin{equation*}
-\frac{b+a}{ab}=v.t.\overset{b}{\underset{-a}{\int }}\frac{dx}{x^{2}}%
.\blacktriangledown
\end{equation*}
\end{example}

Let the singular domain $\vec{r}_{g}$ of the function $f\left( \vec{r}%
\right) $\textit{:} $f\left( \vec{r}\right) =$ $f\left( z^{*},z\right) $,
defined on the complex plane $\vec{\varrho}$, be a cylindrical domain $\vec{r%
}_{\Sigma }$ bounded by contour surface $\vec{r}_{\gamma }$ of the vector
space $\vec{r}$\textit{:} $\vec{r}=\vec{\varrho}+\varkappa \vec{n}$, and
whose bases are obtained by a translation of the domain $\vec{\varrho}_{g}$
bounded by contour $\vec{\varrho}_{\gamma }$ of the complex plane $\vec{%
\varrho}$ to the direction of the unit normal vector $\vec{n}$ for the
constant values\textit{:} $-h$ and $h$\ ($\varkappa _{1}\left( \vec{\varrho}%
\right) =-h$ and $\varkappa _{2}\left( \vec{\varrho}\right) =h$). In this
case, if the function $f\left( \vec{r}\right) $ is integrable over the
contour of integration $\vec{\varrho}_{\gamma }$, then on the basis of the
result (\ref{31}) it follows that 
\begin{equation*}
\overset{\circlearrowleft }{\underset{\vec{\varrho}_{\gamma }}{\int }}%
f\left( z,z^{*}\right) dz-v.p.\underset{\vec{\varrho}_{g}}{\iint }\frac{%
\partial }{\partial z^{*}}f\left( z,z^{*}\right) dz^{*}dz=\underset{\vec{%
\varrho}_{N}\in v.s.\vec{\varrho}_{g}}{\sum }p_{\vec{\varrho}_{\gamma
}\leftrightarrow \vec{\varrho}_{N}}\underset{\vec{\varrho}=\vec{\varrho}_{N}%
}{Res}f\left( z,z^{*}\right) ,
\end{equation*}
\begin{equation*}
-\overset{\circlearrowleft }{\underset{\vec{\varrho}_{\gamma }}{\int }}%
f\left( z,z^{*}\right) dz^{*}-v.p.\underset{\vec{\varrho}_{g}}{\iint }\frac{%
\partial }{\partial z}f\left( z,z^{*}\right) dz^{*}dz=\underset{\vec{\varrho}%
_{N}\in v.s.\vec{\varrho}_{g}}{\sum }p_{\vec{\varrho}_{\gamma
}\leftrightarrow \vec{\varrho}_{N}}\underset{\vec{\varrho}=\vec{\varrho}_{N}%
}{Res^{\star }}f\left( z,z^{*}\right) .
\end{equation*}

Clearly, the partial residues of the function $f\left( z,z^{*}\right) $ are
by definition 
\begin{equation}
\overset{\circlearrowleft }{\underset{\vec{\varrho}_{s}\left( \vec{\varrho}%
_{N},d\delta \right) }{\int }}f\left( z,z^{*}\right) dz=p_{\vec{\varrho}%
_{s}\left( \vec{\varrho}_{N},d\delta \right) \leftrightarrow \vec{\varrho}%
_{N}}\underset{\vec{\varrho}=\vec{\varrho}_{N}}{Res}f\left( z,z^{*}\right) ,
\label{34}
\end{equation}
\begin{equation}
-\overset{\circlearrowleft }{\underset{\vec{\varrho}_{s}\left( \vec{\varrho}%
_{N},d\delta \right) }{\int }}f\left( z,z^{*}\right) dz^{*}=p_{\vec{\varrho}%
_{s}\left( \vec{\varrho}_{N},d\delta \right) \leftrightarrow \vec{\varrho}%
_{N}}\underset{\vec{\varrho}=\vec{\varrho}_{N}}{Res^{\star }}f\left(
z,z^{*}\right) .  \label{35}
\end{equation}

\begin{definition}
\textit{If and only if the function }$f\left( \vec{\varrho}\right) $\textit{%
: }$f\left( \vec{\varrho}\right) =$\textit{\ }$f\left( z,z^{*}\right) $%
\textit{, at each point }$\vec{\varrho}_{N}$\textit{\ of domain }$\vec{%
\varrho}_{g}$\textit{\ of the complex plane }$\vec{\varrho}$\textit{, which
is a regular domain of the function, satisfies one of the conditions: }$%
\left\{ \frac{\partial }{\partial z^{*}}f\left( z,z^{*}\right) \right\} _{%
\vec{\varrho}_{N}}=0$\textit{\ or }$\left\{ \frac{\partial }{\partial z}%
f\left( z,z^{*}\right) \right\} _{\vec{\varrho}_{N}}=0$\textit{, the
function }$f\left( \vec{r}\right) $\textit{\ is a regular-analytic function
on the domain }$\vec{\varrho}_{g}$\textit{.}$\blacktriangledown $
\end{definition}

\begin{definition}
\label{s-a.}\textit{If and only if the function }$f\left( \vec{\varrho}%
\right) $\textit{: }$f\left( \vec{\varrho}\right) =$\textit{\ }$f\left(
z,z^{*}\right) $\textit{, at each point }$\vec{\varrho}_{N}$\textit{\ of
domain }$\vec{\varrho}_{g}$\textit{\ of the complex plane }$\vec{\varrho}$%
\textit{, which is a singular domain of the function, satisfies one of the
conditions: }$\left\{ \frac{\partial }{\partial z^{*}}f\left( z,z^{*}\right)
\right\} _{\vec{\varrho}_{N}}=0$\textit{\ or }$\left\{ \frac{\partial }{%
\partial z}f\left( z,z^{*}\right) \right\} _{\vec{\varrho}_{N}}=0$\textit{,
the function }$f\left( \vec{r}\right) $\textit{\ is a singular-analytic
function on the domain }$\vec{\varrho}_{g}$\textit{.}$\blacktriangledown $
\end{definition}

For a complex vector valued function $\vec{F}\left( \vec{\varrho}\right) $%
\textit{:} $\vec{F}\left( \vec{\varrho}\right) =P\left( \vec{\varrho}\right) 
\vec{w}_{1}+Q\left( \vec{\varrho}\right) \vec{w}_{2}$, whose vector field is
on two-dimensional vector space $\vec{\varrho}$, and functions\textit{:} $%
P\left( \vec{\varrho}\right) $ and $Q\left( \vec{\varrho}\right) $, are
integrable over contour of an integration $\vec{\varrho}_{\gamma }$ bounding
a singular domain $\vec{\varrho}_{g}$ of the function $\vec{F}\left( \vec{%
\varrho}\right) $%
\begin{equation*}
\overset{\circlearrowleft }{\underset{\vec{\varrho}_{\gamma }}{\int }}\left[ 
\vec{F}\left( \vec{\varrho}\right) \times d\vec{\varrho}\right] \cdot \vec{n}%
-v.p.\underset{\vec{\varrho}_{g}}{\iint }\left[ \nabla \cdot \vec{F}\left( 
\vec{\varrho}\right) \right] \left( d\vec{\sigma}\cdot \vec{n}\right) =%
\underset{\vec{\varrho}_{N}\in v.s.\vec{\varrho}_{g}}{\sum }p_{\vec{\varrho}%
_{\gamma }\leftrightarrow \vec{\varrho}_{N}}\underset{\vec{\varrho}=\vec{%
\varrho}_{N}}{Res}\vec{F}\left( \vec{\varrho}\right) ,
\end{equation*}
\begin{equation*}
\overset{\circlearrowleft }{\underset{\vec{\varrho}_{\gamma }}{\int }}\vec{F}%
\left( \vec{\varrho}\right) \cdot d\vec{\varrho}-v.p.\underset{\vec{\varrho}%
_{g}}{\iint }\vec{n}\cdot \left[ \nabla \times \vec{F}\left( \vec{\varrho}%
\right) \right] \left( d\vec{\sigma}\cdot \vec{n}\right) =\underset{\vec{%
\varrho}_{N}\in v.s.\vec{\varrho}_{g}}{\sum }p_{\vec{\varrho}_{\gamma
}\leftrightarrow \vec{\varrho}_{N}}\underset{\vec{\varrho}=\vec{\varrho}_{N}%
}{Res^{\star }}\vec{F}\left( \vec{\varrho}\right) .
\end{equation*}

In this case 
\begin{equation*}
\underset{\vec{\varrho}=\vec{\varrho}_{N}}{Res}\vec{F}\left( \vec{\varrho}%
\right) =\underset{\vec{\varrho}=\vec{\varrho}_{N}}{Res}P\left(
z,z^{*}\right) +\underset{\vec{\varrho}=\vec{\varrho}_{N}}{Res^{\star }}%
Q\left( z,z^{*}\right) ,
\end{equation*}
\begin{equation*}
\underset{\vec{\varrho}=\vec{\varrho}_{N}}{Res^{\star }}\vec{F}\left( \vec{%
\varrho}\right) =\underset{\vec{\varrho}=\vec{\varrho}_{N}}{Res}Q\left(
z,z^{*}\right) -\underset{\vec{\varrho}=\vec{\varrho}_{N}}{Res^{\star }}%
P\left( z,z^{*}\right) .
\end{equation*}

\begin{example}
The domain $\vec{\varrho}_{g_{\delta }}$\textit{:} $\vec{\varrho}_{g_{\delta
}}=\left\{ \vec{\varrho}\text{\textit{:} }a\geq \left| \vec{\varrho}\right|
\geq \delta \text{\textit{;} }\left( \delta ,a\right) \in R_{+}^{1}\right\} $%
, of the complex plane $\vec{\varrho}$, is a regular domain of the function%
\textit{:} $\left( z,z^{*}\right) \mapsto \frac{1}{2}\log $ $\left(
zz^{*}\right) $. If one considers the fact that $\frac{1}{2}\log $ $\left(
zz^{*}\right) =\frac{1}{2}\ln \left( x^{2}+y^{2}\right) $, then on the basis
of the result of well-known \textit{Green-Riemann's} theorem 
\begin{equation*}
\overset{\circlearrowleft }{\underset{\vec{\varrho}_{a}}{\int }}\log \left(
zz^{*}\right) dz-\overset{\circlearrowleft }{\underset{\vec{\varrho}_{\delta
}}{\int }}\log \left( zz^{*}\right) dz=2i\underset{\vec{\varrho}_{g_{\delta
}}}{\iint }\frac{x+iy}{x^{2}+y^{2}}dxdy,
\end{equation*}
\begin{equation*}
-\overset{\circlearrowleft }{\underset{\vec{\varrho}_{a}}{\int }}\log \left(
zz^{*}\right) dz^{*}+\overset{\circlearrowleft }{\underset{\vec{\varrho}%
_{\delta }}{\int }}\log \left( zz^{*}\right) dz^{*}=2i\underset{\vec{\varrho}%
_{g_{\delta }}}{\iint }\frac{x-iy}{x^{2}+y^{2}}dxdy.
\end{equation*}

Similarly, if one takes the fact that $\frac{1}{2}\log $ $\frac{z}{z^{*}}%
=i\arctan \frac{y}{x}$ into account 
\begin{equation*}
\overset{\curvearrowleft }{\underset{\vec{\varrho}_{a}}{\overset{\vec{\varrho%
}_{B}\vec{\varrho}_{A}}{\int }}}\log \frac{z}{z^{*}}dz+\overset{\vec{\varrho}%
_{B}}{\overset{\searrow }{\overset{\vec{\varrho}_{C}}{\underset{\vec{\varrho}%
_{1}}{\int }}}}\log \frac{z}{z^{*}}dz+\overset{\curvearrowright }{\underset{%
\vec{\varrho}_{\delta }}{\overset{\vec{\varrho}_{C}\vec{\varrho}_{D}}{\int }}%
}\log \frac{z}{z^{*}}dz+
\end{equation*}
\begin{equation*}
+\overset{\vec{\varrho}_{D}}{\overset{\swarrow }{\overset{\vec{\varrho}_{A}}{%
\underset{\vec{\varrho}_{2}}{\int }}}}\log \frac{z}{z^{*}}dz=-2i\underset{%
\vec{\varrho}_{g_{\delta }}\backslash \vartriangleright \vec{\varrho}%
_{g_{\delta }}}{\iint }\frac{x+iy}{x^{2}+y^{2}}dxdy,
\end{equation*}
\begin{equation*}
\overset{\curvearrowleft }{\underset{\vec{\varrho}_{a}}{\overset{\vec{\varrho%
}_{B}\vec{\varrho}_{A}}{\int }}}\log \frac{z}{z^{*}}dz^{*}+\overset{\vec{%
\varrho}_{B}}{\overset{\searrow }{\overset{\vec{\varrho}_{C}}{\underset{\vec{%
\varrho}_{1}}{\int }}}}\log \frac{z}{z^{*}}dz^{*}+\overset{\curvearrowright 
}{\underset{\vec{\varrho}_{\delta }}{\overset{\vec{\varrho}_{C}\vec{\varrho}%
_{D}}{\int }}}\log \frac{z}{z^{*}}dz^{*}+
\end{equation*}
\begin{equation*}
+\overset{\vec{\varrho}_{D}}{\overset{\swarrow }{\overset{\vec{\varrho}_{A}}{%
\underset{\vec{\varrho}_{2}}{\int }}}}\log \frac{z}{z^{*}}dz^{*}=-2i%
\underset{\vec{\varrho}_{g_{\delta }}\backslash \vartriangleright \vec{%
\varrho}_{g_{\delta }}}{\iint }\frac{x-iy}{x^{2}+y^{2}}dxdy,
\end{equation*}
where the domain $\vec{\varrho}_{g_{\delta }}\backslash \vartriangleright 
\vec{\varrho}_{g_{\delta }}$ is a part of the domain $\vec{\varrho}%
_{g_{\delta }}$ bounded by parts of circular contours of an integration%
\textit{:} $\vec{\varrho}_{a}$ and $\vec{\varrho}_{\delta }$, bounding the
domain $\vec{\varrho}_{g_{\delta }}$ and by the segments of straight-lines $%
\vec{\varrho}_{k}$ of the complex plane $\vec{\varrho}$\textit{: }$\vec{%
\varrho}_{k}=\left\{ \vec{\varrho}\text{\textit{:} }\,\vec{\varrho}=\left| 
\vec{\varrho}\right| \vec{\varrho}_{ok}\text{ }(\varphi =\varphi
_{k})\right\} $, ($k=1,2$). The points\textit{:} $\vec{\varrho}_{A}$,$\vec{%
\varrho}_{B}$, $\vec{\varrho}_{C}$ and $\vec{\varrho}_{D}$, are points
obtained by an intersection of circular contours of an integration\textit{:} 
$\vec{\varrho}_{a}$ and $\vec{\varrho}_{\delta }$, with directions $\vec{%
\varrho}_{k}$.

For arbitrary chosen angular values $\varphi _{k}$, and in the case as%
\textit{:} $\varphi _{1}\rightarrow \pi $ and $\varphi _{2}\rightarrow -\pi $%
, 
\begin{equation*}
\overset{\circlearrowleft }{\underset{\vec{\varrho}_{a}}{\int }}\log \frac{z%
}{z^{*}}dz-\overset{\circlearrowleft }{\underset{\vec{\varrho}_{\delta }}{%
\int }}\log \frac{z}{z^{*}}dz+\overset{\rightleftharpoons }{\underset{\vec{%
\varrho}_{k}}{\int }}\log \frac{z}{z^{*}}dz=-2i\underset{\vec{\varrho}%
_{g_{\delta }}}{\iint }\frac{x+iy}{x^{2}+y^{2}}dxdy,
\end{equation*}
\begin{equation*}
\overset{\circlearrowleft }{\underset{\vec{\varrho}_{a}}{\int }}\log \frac{z%
}{z^{*}}dz^{*}-\overset{\circlearrowleft }{\underset{\vec{\varrho}_{\delta }%
}{\int }}\log \frac{z}{z^{*}}dz^{*}+\overset{\rightleftharpoons }{\underset{%
\vec{\varrho}_{k}}{\int }}\log \frac{z}{z^{*}}dz^{*}=-2i\underset{\vec{%
\varrho}_{g_{\delta }}}{\iint }\frac{x-iy}{x^{2}+y^{2}}dxdy.
\end{equation*}

In other words the scalar valued function $z\mapsto \log z$\textit{:} 
\begin{equation*}
\log z=\frac{1}{2}\left[ \log \left( zz^{*}\right) +\log \frac{z}{z^{*}}%
\right] ,
\end{equation*}
is a singular-analytic function on the domain $\vec{\varrho}_{g}$\textit{:} $%
\vec{\varrho}_{g}=\left\{ \vec{\varrho}\text{\textit{:} }\left| \vec{\varrho}%
\right| \leq a\right\} $, of complex plane $\vec{\varrho}$, more exactly 
\begin{equation*}
\overset{\circlearrowleft }{\underset{\vec{\varrho}_{a}}{\int }}\log zdz+%
\overset{\rightleftharpoons }{\underset{\vec{\varrho}_{k}}{\int }}\log zdz=%
\underset{\delta \rightarrow 0^{+}}{\lim }\overset{\circlearrowleft }{%
\underset{\vec{\varrho}_{\delta }}{\int }}\log zdz,
\end{equation*}
\begin{equation*}
-\overset{\circlearrowleft }{\underset{\vec{\varrho}_{a}}{\int }}\log
zdz^{*}-v.p.\underset{\vec{\varrho}_{g}}{\iint }\frac{1}{z}dzdz^{*}-\overset{%
\rightleftharpoons }{\underset{\vec{\varrho}_{k}}{\int }}\log zdz^{*}=-%
\underset{\delta \rightarrow 0^{+}}{\lim }\overset{\circlearrowleft }{%
\underset{\vec{\varrho}_{\delta }}{\int }}\log zdz^{*}.\blacktriangledown
\end{equation*}
\end{example}

\section{The fundamental lemmas}

Let the complex plane $\vec{\varrho}$ be subdivided by finitely many
directions $\vec{\varrho}_{k}$\textit{:} 
\begin{equation*}
\vec{\varrho}_{k}=\left\{ \vec{\varrho}\text{\textit{:} }\,\vec{\varrho}-%
\vec{\varrho}_{N}=\left| \vec{\varrho}-\vec{\varrho}_{N}\right| \vec{\varrho}%
_{ok}^{\vec{\varrho}_{N}}\text{ }(\varphi =\varphi _{k})\right\} ,
\end{equation*}
where $\vec{\varrho}_{ok}^{\vec{\varrho}_{N}}$ is an unit vector $\vec{%
\varrho}_{ok}$ at the point $\vec{\varrho}_{N}$ of the complex plane $\vec{%
\varrho}$, into $K$ ($k=1,2,...,K$) different domains of convergence $\vec{%
\varrho}_{g_{k}}$\textit{:} 
\begin{equation*}
\vec{\varrho}_{g_{k}}=\left\{ \vec{\varrho}\text{\textit{:} }\,\vec{\varrho}-%
\vec{\varrho}_{N}=\left| \vec{\varrho}-\vec{\varrho}_{N}\right| \vec{\varrho}%
_{o}^{\vec{\varrho}_{N}},\varphi _{k+1}>\varphi >\varphi _{k}\text{ }%
(\varphi _{K+1}=2\pi )\right\} ,
\end{equation*}
of the function $\vec{\varrho}\mapsto $ $\left( \vec{\varrho}-2z^{*}\vec{w}%
_{1}\right) f\left( \vec{\varrho}+\vec{\varrho}_{N}\right) $\textit{: } 
\begin{equation}
\underset{\left| \vec{\varrho}\right| \rightarrow 0^{+}}{\lim }\left( \vec{%
\varrho}-2z^{*}\vec{w}_{1}\right) f\left( \vec{\varrho}+\vec{\varrho}%
_{N}\right) =\vec{A}_{0k};\text{ }\vec{\varrho}\in \vec{\varrho}_{g_{k}}.
\label{36}
\end{equation}

In that case $\delta \left( \vec{\varrho}_{N}\right) $-neighborhood of the
point $\vec{\varrho}_{N}$, bounded by a circular path of an
integration\thinspace $\vec{\varrho}_{\delta }$ centred at the point $\vec{%
\varrho}_{N}$ and of an arbitrary small radius $\delta \left( \vec{\varrho}%
_{N}\right) $, and over which the function is integrable by assumption, is
subdivided by the direction $\vec{\varrho}_{k}$ into sub-domains $\vec{%
\varrho}_{g_{k}}^{\delta }$\textit{:} 
\begin{equation*}
\vec{\varrho}_{g_{k}}^{\delta }=\left\{ \vec{\varrho}\text{\textit{:} }\vec{%
\varrho}-\vec{\varrho}_{N}\leq \delta \vec{\varrho}_{o}^{\vec{\varrho}%
_{N}},\varphi _{k+1}>\varphi >\varphi _{k}\right\} .
\end{equation*}

For every sub-domain $\vec{\varrho}_{g_{k}}^{\delta \varepsilon }$ of the
domain of convergence $\vec{\varrho}_{g_{k}}$, bounded by contour $\vec{%
\varrho}_{\gamma _{k}}^{\delta \varepsilon }$\textit{:} 
\begin{equation*}
\vec{\varrho}_{\gamma _{k}}^{\delta \varepsilon }=\vec{\varrho}%
_{k+\varepsilon }\overset{\vec{\varrho}_{\delta }}{\underset{\vec{\varrho}%
_{N}}{\nearrow }}\cup \overset{\curvearrowleft }{\vec{\varrho}_{\delta _{k}}}%
\cup \vec{\varrho}_{\left( k+1\right) -\varepsilon }\overset{\vec{\varrho}%
_{\delta }}{\underset{\vec{\varrho}_{N}}{\searrow }},
\end{equation*}
where 
\begin{equation*}
\vec{\varrho}_{k+\varepsilon }=\left\{ \vec{\varrho}\text{\textit{:} }\sqrt{2%
}\,\left( \vec{\varrho}-\vec{\varrho}_{N}\right) =\left| \vec{\varrho}-\vec{%
\varrho}_{N}\right| \left\{ e^{-i\left[ \varphi _{k}+\varepsilon \left(
\delta \right) \right] }\vec{w}_{1}+e^{i\left[ \varphi _{k}+\varepsilon
\left( \delta \right) \right] }\vec{w}_{2}\right\} ^{\vec{\varrho}%
_{N}}\right\} ,
\end{equation*}
\begin{equation*}
\vec{\varrho}_{\left( k+1\right) -\varepsilon }=\left\{ \vec{\varrho}\text{%
\textit{:} }\,\sqrt{2}\left( \vec{\varrho}-\vec{\varrho}_{N}\right) =\left| 
\vec{\varrho}-\vec{\varrho}_{N}\right| \left\{ e^{-i\left[ \varphi
_{k+1}-\varepsilon \left( \delta \right) \right] }\vec{w}_{1}+e^{i\left[
\varphi _{k+1}-\varepsilon \left( \delta \right) \right] }\vec{w}%
_{2}\right\} ^{\vec{\varrho}_{N}}\right\} ,
\end{equation*}
and the function $\varepsilon \left( \delta \right) \in R_{+}^{1}$, being
finite small values, satisfies the condition\textit{: }$\underset{\delta
\rightarrow 0^{+}}{\lim }\varepsilon \left( \delta \right) =0$, it holds 
\begin{equation*}
\mathit{\,}\underset{\delta \rightarrow 0^{+}}{\lim }\overset{%
\curvearrowleft }{\underset{\vec{\varrho}_{\delta _{k}}}{\int }}f\left( \vec{%
\varrho}\right) d\vec{\varrho}=i\underset{\delta \rightarrow 0^{+}}{\lim }%
\overset{\varphi _{k+1}-\varepsilon \left( \delta \right) }{\underset{%
\varphi _{k}+\varepsilon \left( \delta \right) }{\int }}\left( \vec{\varrho}%
-2z^{*}\vec{w}_{1}\right) f\left( \vec{\varrho}+\vec{\varrho}_{N}\right)
d\varphi =\alpha _{k}i\vec{A}_{0k},
\end{equation*}
where $\alpha _{k}=\varphi _{k+1}-\varphi _{k}$. Accordingly, and on the
basis of equality (\ref{26}) of \textit{Definition \ref{res inf.},} finally
it follows that 
\begin{equation}
2\pi i\underset{\vec{\varrho}=\vec{\varrho}_{N}}{\overset{\rightarrow }{Res}}%
f\left( \vec{\varrho}\right) =\overset{\circlearrowleft }{\underset{\vec{%
\varrho}_{s}\left( \vec{\varrho}_{N},d\delta \right) }{\int }}f\left( \vec{%
\varrho}\right) d\vec{\varrho}=\overset{K}{\underset{k=1}{\sum }}\alpha _{k}i%
\vec{A}_{0k},  \label{37}
\end{equation}
where $\underset{\vec{\varrho}=\vec{\varrho}_{N}}{\overset{\rightarrow }{Res}%
}f\left( \vec{\varrho}\right) =-\underset{\vec{\varrho}=\vec{\varrho}_{N}}{%
Res}^{\star }f\left( \vec{\varrho}\right) \vec{w}_{1}+$ $\underset{\vec{%
\varrho}=\vec{\varrho}_{N}}{Res}f\left( \vec{\varrho}\right) \vec{w}_{2}$.

On the other hand, on the basis of the integral equality 
\begin{equation*}
\left( \vec{w}_{1}+\vec{w}_{2}\right) \cdot \overset{\circlearrowleft }{%
\underset{\vec{\varrho}_{\delta }}{\int }}f\left( \vec{\varrho}\right) d\vec{%
\varrho}=\overset{\circlearrowleft }{\underset{\vec{\varrho}_{\zeta }}{\int }%
}\vec{\varrho}^{-2}f\left( \vec{\varrho}^{-1}+\vec{\varrho}_{N}\right) \cdot
d\vec{\varrho},
\end{equation*}
where\textit{:} $\vec{\varrho}^{-1}=\left( z^{*}\right) ^{-1}\vec{w}%
_{1}+z^{-1}\vec{w}_{2}$ and $\vec{\varrho}^{-2}=\left( z^{*}\right) ^{-2}%
\vec{w}_{1}+z^{-2}\vec{w}_{2}$ as well as $\vec{\varrho}_{\zeta }=\left\{ 
\vec{\varrho}\text{\textit{:} }\vec{\varrho}=\frac{1}{\delta }\vec{\varrho}%
_{o}\right\} $, more exactly on the basis of the integral equality 
\begin{equation*}
\overset{2\pi }{\underset{0}{\int }}\left( \vec{\varrho}-2z^{*}\vec{w}%
_{1}\right) f\left( \vec{\varrho}+\vec{\varrho}_{N}\right) d\varphi =%
\overset{2\pi }{\underset{0}{\int }}\left( \vec{\varrho}^{-1}-\frac{2}{z^{*}}%
\vec{w}_{1}\right) f\left( \vec{\varrho}^{-1}+\vec{\varrho}_{N}\right)
d\varphi ,
\end{equation*}
it follows that 
\begin{equation*}
\left( \vec{w}_{1}+\vec{w}_{2}\right) \cdot \underset{\delta \rightarrow
0^{+}}{\lim }\overset{2\pi }{\underset{0}{\int }}\left( \vec{\varrho}^{-1}-%
\frac{2}{z^{*}}\vec{w}_{1}\right) f\left( \vec{\varrho}^{-1}+\vec{\varrho}%
_{N}\right) d\varphi =
\end{equation*}
\begin{equation*}
=2\pi i\underset{\vec{\varrho}_{M}\in \vec{\varrho}}{\sum }\underset{\vec{%
\varrho}=\vec{\varrho}_{M}}{Res^{*}}\left[ \vec{\varrho}^{-2}f\left( \vec{%
\varrho}^{-1}+\vec{\varrho}_{N}\right) \right] .
\end{equation*}

In other words, if the condition (\ref{35}) is satisfied, more exactly the
condition being an equivalent to it 
\begin{equation}
\underset{\left| \vec{\varrho}\right| \rightarrow +\infty }{\lim }\left( 
\vec{\varrho}^{-1}-\frac{2}{z^{*}}\vec{w}_{1}\right) f\left( \vec{\varrho}%
^{-1}+\vec{\varrho}_{N}\right) =\vec{A}_{\infty k};\vec{\varrho}\in \vec{%
\varrho}_{g_{k}},  \label{38}
\end{equation}
where $\vec{A}_{\infty k}=\vec{A}_{0k}$, then 
\begin{equation}
-2\pi i\underset{\vec{\varrho}=\vec{\varrho}_{\infty }}{Res^{*}}\left[ \vec{%
\varrho}^{-2}f\left( \vec{\varrho}^{-1}+\vec{\varrho}_{N}\right) \right]
=\left( \vec{w}_{1}+\vec{w}_{2}\right) \cdot \overset{K}{\underset{k=1}{\sum 
}}\alpha _{k}i\vec{A}_{\infty k},  \label{39}
\end{equation}
that is 
\begin{equation*}
-\underset{\vec{\varrho}=\vec{\varrho}_{\infty }}{Res^{*}}\left[ \vec{\varrho%
}^{-2}f\left( \vec{\varrho}^{-1}+\vec{\varrho}_{N}\right) \right] =\left( 
\vec{w}_{1}+\vec{w}_{2}\right) \cdot \underset{\vec{\varrho}=\vec{\varrho}%
_{N}}{\overset{\rightarrow }{Res}}f\left( \vec{\varrho}\right) .
\end{equation*}

Clearly, in the inverse case, if the condition is satisfied 
\begin{equation}
\underset{\left| \vec{\varrho}\right| \rightarrow +\infty }{\lim }\left( 
\vec{\varrho}-2z^{*}\vec{w}_{1}\right) f\left( \vec{\varrho}+\vec{\varrho}%
_{N}\right) =\vec{A}_{\infty k};\vec{\varrho}\in \vec{\varrho}_{g_{k}},
\label{40}
\end{equation}
then 
\begin{equation}
\underset{\vec{\varrho}=\vec{\varrho}_{\infty }}{\overset{\rightarrow }{Res}}%
f\left( \vec{\varrho}\right) =-\overset{K}{\underset{k=1}{\sum }}\alpha _{k}i%
\vec{A}_{\infty k},
\end{equation}
more exactly 
\begin{equation}
-\left( \vec{w}_{1}+\vec{w}_{2}\right) \cdot \underset{\vec{\varrho}=\vec{%
\varrho}_{\infty }}{\overset{\rightarrow }{Res}}f\left( \vec{\varrho}\right)
=\underset{\vec{\varrho}=\vec{\varrho}_{0}}{Res^{*}}\left[ \vec{\varrho}%
^{-2}f\left( \vec{\varrho}^{-1}+\vec{\varrho}_{N}\right) \right] .
\label{42}
\end{equation}

\begin{lemma}
\label{3.1}\textit{Let directions }$\vec{\varrho}_{k}$\textit{: }$\vec{%
\varrho}_{k}=\left\{ \vec{\varrho}\text{\textit{:} }\,\vec{\varrho}-\vec{%
\varrho}_{C}=\left| \vec{\varrho}-\vec{\varrho}_{C}\right| \vec{\varrho}%
_{ok}^{\vec{\varrho}_{C}}\right\} $\textit{\ (}$k=1,2,...,K$\textit{)
subdivide the complex plane }$\vec{\varrho}$\textit{\ into }$K$\textit{: }$%
K=2$\textit{, domain }$\vec{\varrho}_{g_{k}}$\textit{: } 
\begin{equation*}
\vec{\varrho}_{g_{k}}=\left\{ \vec{\varrho}\text{\textit{:} }\,\vec{\varrho}-%
\vec{\varrho}_{C}=\left| \vec{\varrho}-\vec{\varrho}_{C}\right| \vec{\varrho}%
_{o}^{\vec{\varrho}_{C}},\varphi _{k+1}>\varphi >\varphi _{k}\text{ }%
(\varphi _{K+1}=2\pi )\right\}
\end{equation*}
\textit{and the function }$\vec{\varrho}\mapsto f\left( \vec{\varrho}\right) 
$\textit{\ satisfies the condition } 
\begin{equation*}
\underset{\left| \vec{\varrho}\right| \rightarrow +\infty }{\lim }\left( 
\vec{\varrho}-2z^{*}\vec{w}_{1}\right) f\left( \vec{\varrho}+\vec{\varrho}%
_{C}\right) =\vec{A}_{\infty 1},\mathit{\ }\vec{\varrho}\in \vec{\varrho}%
_{g_{1}}.
\end{equation*}

\textit{In that case} 
\begin{equation}
\overset{\curvearrowleft }{\underset{ext.\vec{\varrho}_{\infty }}{\int }}%
f\left( \vec{\varrho}\right) d\vec{\varrho}=2\pi i\underset{\vec{\varrho}%
_{N}\in \vec{\varrho}}{\sum }\underset{\vec{\varrho}=\vec{\varrho}_{N}}{%
\overset{\rightarrow }{Res}}f\left( \vec{\varrho}\right) -\alpha _{1}i\vec{A}%
_{\infty 1},  \label{43}
\end{equation}
\textit{where }$ext.\vec{\varrho}_{\infty }$\textit{\ is a set of an
infinite points outside of the domain }$\vec{\varrho}_{g_{1}}$\textit{, and
angle }$\alpha _{1}$\textit{\ is equal to an angular difference: }$\alpha
_{1}=\varphi _{2}-\varphi _{1}$\textit{.}$\blacktriangledown $
\end{lemma}

Preceding lemma is an explicit consequence of the equality (\ref{26}) of 
\textit{Definition \ref{res inf.}} and of the result (\ref{41}). On the
basis of the result (\ref{37}) it can be formulated the lemma being
analogous to \textit{Lemma \ref{3.1}}

\begin{lemma}
\label{3.2}\textit{Let directions }$\vec{\varrho}_{k}$\textit{: }$\vec{%
\varrho}_{k}=\left\{ \vec{\varrho}\text{\textit{:} }\,\vec{\varrho}-\vec{%
\varrho}_{N}=\left| \vec{\varrho}-\vec{\varrho}_{N}\right| \vec{\varrho}%
_{ok}^{\vec{\varrho}_{N}}\right\} $ ($k=1,2,...,K$) \textit{subdivide the
complex plane }$\vec{\varrho}$\textit{\ into }$K$\textit{: }$K=2$\textit{,
domain }$\vec{\varrho}_{g_{k}}$\textit{: } 
\begin{equation*}
\vec{\varrho}_{g_{k}}=\left\{ \vec{\varrho}\text{\textit{:} }\,\vec{\varrho}-%
\vec{\varrho}_{N}=\left| \vec{\varrho}-\vec{\varrho}_{N}\right| \vec{\varrho}%
_{o}^{\vec{\varrho}_{N}},\varphi _{k+1}>\varphi >\varphi _{k}\text{ }%
(\varphi _{K+1}=2\pi )\right\}
\end{equation*}
\textit{and the function }$\vec{\varrho}\mapsto f\left( \vec{\varrho}\right) 
$\textit{\ satisfies the condition } 
\begin{equation*}
\underset{\left| \vec{\varrho}\right| \rightarrow 0^{+}}{\lim }\left( \vec{%
\varrho}-2z^{*}\vec{w}_{1}\right) f\left( \vec{\varrho}+\vec{\varrho}%
_{N}\right) =\vec{A}_{01},\mathit{\ }\vec{\varrho}\in \vec{\varrho}_{g_{1}}.
\end{equation*}

\textit{In that case} 
\begin{equation}
\overset{\curvearrowleft }{\underset{ext.\vec{\varrho}_{s}\left( \vec{\varrho%
}_{N},d\delta \right) }{\int }}f\left( \vec{\varrho}\right) d\vec{\varrho}%
=2\pi i\underset{\vec{\varrho}=\vec{\varrho}_{N}}{\overset{\rightarrow }{Res}%
}f\left( \vec{\varrho}\right) -\alpha _{1}i\vec{A}_{01},  \label{44}
\end{equation}
\textit{where }$ext.\vec{\varrho}_{s}\left( \vec{\varrho}_{N},d\delta
\right) $\textit{\ is a part of infinitesimally small circular path of an
integration }$\vec{\varrho}_{s}\left( \vec{\varrho}_{N},d\delta \right) $%
\textit{\ outside of the domain }$\vec{\varrho}_{g_{1}}$\textit{, and angle }%
$\alpha _{1}$\textit{\ is equal to an angular difference: }$\alpha
_{1}=\varphi _{2}-\varphi _{1}$\textit{.}$\blacktriangledown $
\end{lemma}

As a consequence either of the result of\textit{\ Lemma \ref{3.1} }or of the
result (\ref{33}), the following result is obtained

\begin{lemma}
\label{3.3}\textit{Let directions }$\vec{\varrho}_{k}$\textit{: }$\vec{%
\varrho}_{k}=\left\{ \vec{\varrho}\text{\textit{:} }\,\vec{\varrho}-\vec{%
\varrho}_{C}=\left| \vec{\varrho}-\vec{\varrho}_{C}\right| \vec{\varrho}%
_{ok}^{\vec{\varrho}_{C}}\right\} $\textit{\ }($k=1,2,...,K$) \textit{%
subdivide the complex plane }$\vec{\varrho}$\textit{\ into }$K$\textit{: }$%
K=2$\textit{, domain }$\vec{\varrho}_{g_{k}}$\textit{: } 
\begin{equation*}
\vec{\varrho}_{g_{k}}=\left\{ \vec{\varrho}\text{\textit{:} }\,\vec{\varrho}-%
\vec{\varrho}_{C}=\left| \vec{\varrho}-\vec{\varrho}_{C}\right| \vec{\varrho}%
_{o}^{\vec{\varrho}_{C}},\varphi _{k+1}>\varphi >\varphi _{k}\text{ }%
(\varphi _{K+1}=2\pi )\right\}
\end{equation*}
\textit{and the function }$\vec{\varrho}\mapsto f\left( \vec{\varrho}\right) 
$\textit{\ satisfies the condition } 
\begin{equation*}
\underset{\left| \vec{\varrho}\right| \rightarrow +\infty }{\lim }\left( 
\vec{\varrho}-2z^{*}\vec{w}_{1}\right) f\left( \vec{\varrho}+\vec{\varrho}%
_{C}\right) =\vec{0},\mathit{\ }\vec{\varrho}\in \vec{\varrho}_{g2}.
\end{equation*}

\textit{In that case} 
\begin{equation}
v.t.\underset{\overset{\overset{\vec{\varrho}_{2}}{\nwarrow }\overset{\vec{%
\varrho}_{1}}{\swarrow }}{\vec{\varrho}_{C}}}{\int }f\left( \vec{\varrho}%
\right) d\vec{\varrho}-v.p.\underset{\overset{\overset{\vec{\varrho}_{2}}{%
\nwarrow }\vec{\varrho}_{g_{2}}\overset{\vec{\varrho}_{1}}{\swarrow }}{\vec{%
\varrho}_{C}}}{\iint }\left[ \left( \vec{n}\times \nabla \right) f\left( 
\vec{\varrho}\right) \right] \left( d\vec{\sigma}\cdot \vec{n}\right) =
\label{45}
\end{equation}
\begin{equation*}
=\underset{\vec{\varrho}_{N}\in v.s.\overset{\overset{\vec{\varrho}_{2}}{%
\nwarrow }\vec{\varrho}_{g_{2}}\overset{\vec{\varrho}_{1}}{\swarrow }}{\vec{%
\varrho}_{C}}}{\sum }p_{\overset{\overset{\vec{\varrho}_{2}}{\nwarrow }%
\overset{\vec{\varrho}_{1}}{\swarrow }}{\vec{\varrho}_{C}}\leftrightarrow 
\vec{\varrho}_{N}}\underset{\vec{\varrho}=\vec{\varrho}_{N}}{\overset{%
\rightarrow }{Res}}f\left( \vec{\varrho}\right) ,
\end{equation*}
\textit{where }$\overset{\overset{\vec{\varrho}_{2}}{\nwarrow }\vec{\varrho}%
_{g_{2}}\overset{\vec{\varrho}_{1}}{\swarrow }}{\vec{\varrho}_{C}}=\vec{%
\varrho}_{g_{2}}\cup \overset{\overset{\vec{\varrho}_{2}}{\nwarrow }\overset{%
\vec{\varrho}_{1}}{\swarrow }}{\vec{\varrho}_{C}}$.$\blacktriangledown $
\end{lemma}

\section{Conclusion}

Based on the defined concept of an absolute integral sum of a complex
function, which is more general with respect to the concept of an integral
sum of ``ordinary'' integral calculus, the result (\ref{10}) is derived
whose an importance if one takes into consideration the fact that it is
immediately derived from defining equality of a sequence of absolute
integral sums $\vec{A}_{\vartriangle _{j_{1},...,j_{n}}\vec{r}_{\gamma }}$
(from defining equality (\ref{7}) of \textit{Definition \ref{spatial}}) is a
general, more exactly its importance is not conditioned by a convergence of
reduced absolute integral sums. In other words, the result (\ref{10}) is a
more general with respect to the result of the well-known \textit{Cauchy's}
fundamental theorem on the residues of \textit{Cauchy's} calculus of
residues. Accordingly, and on the basis of the redefined concept of a
residue of a complex function as well as of the defined concept of a total
value of an improper integral of the function over a certain singular domain
bounded by a contour surface of the three-dimensional vector space, the
results being more general with respect to the fundamental results of 
\textit{Cauchy's} calculus of residues of both analytic and non-analytic
function on the complex plane, are obtained, the results\textit{:} (\ref{31}%
) and (\ref{33}), as well as the results of the \textit{Section \ref{sec}}.
The concept itself of analytic functions is also redefined, so that in the
class of singular-analytic functions defined by \textit{Definition \ref{s-a.}%
,} in addition to the class of standard analytic functions the class of
pseudo-analytic functions is also contained.

Obtained results give the solid base to a further generalization, whether
they are the fundamentals or not, of the results of \textit{Cauchy's}
calculus of residues on the one hand, as well as of the results in another
areas of both mathematics and applied mathematics, in which the results of 
\textit{Cauchy's} calculus of residues are made use, on the other hand.

\end{document}